\documentclass[12pt,reqno]{amsart}

\usepackage{amsmath}
\usepackage{amsxtra}
\usepackage{amscd}
\usepackage{amsthm}
\usepackage{amsfonts}
\usepackage{amssymb}
\usepackage{eucal}

\newtheorem{theorem}{Theorem}[section]
\newtheorem{conj}[theorem]{Conjecture}
\newtheorem{cor}[theorem]{Corollary}
\newtheorem{lem}[theorem]{Lemma}
\newtheorem{prop}[theorem]{Proposition}

\theoremstyle{definition}

\theoremstyle{remark}

\theoremstyle{remark}

\numberwithin{equation}{section}

\newcommand{\nc}{\newcommand}
\nc{\on}{\operatorname}
\nc{\ch}{\mbox{ch}}
\nc{\Z}{{\mathbb Z}}
\nc{\C}{{\mathbb C}}
\nc{\pone}{{\mathbb C}{\mathbb P}^1}
\nc{\pa}{\partial}
\nc{\F}{{\mathcal F}}
\nc{\arr}{\rightarrow}
\nc{\larr}{\longrightarrow}
\nc{\al}{\alpha}
\nc{\ri}{\rangle}
\nc{\lef}{\langle}
\nc{\W}{{\mathcal W}}
\nc{\la}{\lambda}
\nc{\ep}{\epsilon}

\nc{\su}{\widehat{{\mathfrak sl}}_2}
\nc{\sw}{{\mathfrak s}{\mathfrak l}}

\nc{\g}{{\mathfrak g}}
\nc{\h}{{\mathfrak h}}
\nc{\n}{{\mathfrak n}}
\nc{\N}{\widehat{\n}}
\nc{\G}{\widehat{\g}}
\nc{\De}{\Delta_+}
\nc{\gt}{\widetilde{\g}}
\nc{\Ga}{\Gamma}
\nc{\one}{{\mathbf 1}}
\nc{\z}{{\mathfrak Z}}
\nc{\zz}{{\mathcal Z}}
\nc{\Hh}{{\mathcal H}_\beta}
\nc{\qp}{q^{\frac{k}{2}}}
\nc{\qm}{q^{-\frac{k}{2}}}
\nc{\La}{\Lambda}
\nc{\wt}{\widetilde}
\nc{\qn}{\frac{[m]_q^2}{[2m]_q}}
\nc{\cri}{_{\on{cr}}}
\nc{\kk}{h^\vee}
\nc{\sun}{\widehat{\sw}_N}
\nc{\hh}{\widehat{\mathfrak h}}
\nc{\HH}{{\mathcal H}_{q,t}}
\nc{\ca}{\wt{{\mathcal A}}_{h,k}(\sw_2)}
\nc{\gl}{\widehat{{\mathfrak g}{\mathfrak l}}_2}
\nc{\el}{\ell}
\nc{\s}{{\mathbf s}}
\nc{\bi}{\bibitem}
\nc{\om}{\omega}
\nc{\WW}{\W_\beta}
\nc{\scr}{{\mathbf S}}
\nc{\ab}{{\mathbf a}}
\nc{\rr}{r}
\nc{\ol}{\overline}
\nc{\con}{qt^{-1} + q^{-1}t}
\nc{\den}{q^{\el-1} t^{-\el+1}+ q^{-\el+1} t^{\el-1}}
\nc{\ds}{\displaystyle}
\nc{\B}{B}
\nc{\A}{{\mathbb A}}
\nc{\GG}{{\mathcal G}}
\nc{\UU}{{\mathcal U}}
\nc{\MM}{{\mathcal M}}
\nc{\CC}{{\mathcal C}}
\nc{\GL}{{}^L G}
\nc{\dzz}{\frac{dz}{z}}
\nc{\Res}{\on{Res}}
\nc{\rep}{{\mathcal R}ep \;}
\nc{\uqg}{U_q \G}
\nc{\uqgg}{U_q \g}
\nc{\Fq}{{\mathbb F}_q}

\nc{\stimes}{\ltimes}
\nc{\K}{\hat{\mathcal K}}
\nc{\Ql}{\ol{\mathbb Q}_\ell}

\nc{\ga}{\gamma}
\nc{\PL}{{}^L P}
\nc{\E}{\mc E}
\nc{\mc}{\mathcal}
\nc{\mbf}{\mathbf}
\nc{\bb}{{\mathfrak b}}
\nc{\OO}{{\mc O}}
\nc{\Po}{{\mc P}}
\nc{\V}{{\mc V}}
\nc{\yy}{{\mc Y}}
\nc{\M}{\mathcal M}
\nc{\Coh}{{{\mathcal C}oh}}
\nc{\Cohn}{\Coh_n}
\nc{\f}{{\mathcal F}}
\nc{\si}{_E}
\nc{\Gaf}{{\mathbb G}_{a,\Fq}}
\nc{\KK}{{\mathfrak k}}

\nc{\PCr}{{ \bs P  (\C[x])^r   }}
\nc{\PCN}{{ \bs P  (\C[x])^N   }}

\nc{\sN}{sl_{2N+1}}

\nc{\Pzr}{{ \bs P(\C((x-z)))^r}}

\nc{\PzN}{{ \bs P(\C((x-z)))^N}}

\newcommand{\bean}{\begin{eqnarray}}
\newcommand{\eean}{\end{eqnarray}}
\newcommand{\be}{\begin{displaymath}}
\newcommand{\ee}{\end{displaymath}}
\newcommand{\bea}{\begin{eqnarray*}}   
\newcommand{\eea}{\end{eqnarray*}}
\newcommand{\bs}{\boldsymbol}
\newcommand{\Ref}[1]{{$($\ref{#1}$)$}}

\begin{document}

\title
{Solutions to the XXX type Bethe ansatz equations and flag varieties}
\author
{E. Mukhin  and A. Varchenko}
\thanks{Research of E.M. is supported in part by NSF grant DMS-0140460.
Research of A.V. is supported in part by NSF grant DMS-9801582}
\address{E.M.: Department of Mathematical Sciences, Indiana University -
Purdue University Indianapolis, 402 North Blackford St, Indianapolis,
IN 46202-3216, USA, \newline mukhin@math.iupui.edu}
\address{A.V.: Department of Mathematics, University of North Carolina 
at Chapel Hill, Chapel Hill, NC 27599-3250, USA, anv@email.unc.edu}

\begin{abstract}
We consider a version of the $A_{N}$ Bethe equation of $XXX$ type and introduce
a reproduction procedure constructing new solutions of this equation
from a given one. 
The set of all solutions obtained from a given one is called a population. 
 
We show that a population is isomorphic to the $sl_{N+1}$ flag
variety and that the populations are in one-to-one correspondence 
with intersection points of suitable Schubert cycles in a Grassmanian variety.

We also obtain similar results for the root systems $B_N$ and
$C_N$. Populations of $B_N$ and $C_N$ type are isomorphic to the flag
varieties of $C_N$ and $B_N$ types respectively.
\end{abstract}

\maketitle

\section{Introduction}
In this paper we consider a system of algebraic equations, see below \Ref{Bethe
  slN}, which we consider as an $A_N$ version of the $XXX$ Bethe equation.
In the simplest case of $N=1$, our system does coincide with the famous and
much studied Bethe equation for the inhomogeneous 
$XXX$-model see \cite{BIK}, \cite{Fd}, \cite{FT} and references therein. 

We call solutions of our system \Ref{Bethe slN} with additional simple
conditions $h$-critical points and study the 
problem of counting the number of the $h$-critical points. 

The definition of the critical points depends 
on the complex non-zero step $h$, complex distinct spectral parameters 
$z_1,\dots,z_n$, $sl_N$ dominant
integral
weights $\La_1,\dots,\La_n$ and another $sl_N$ weight $\La_\infty$. 
We conjecture that for generic
positions of $z_i$ and dominant integral $\La_\infty$ 
the number of $h$-critical points modulo obvious symmetries
equals the multiplicity of $L_{\La_\infty}$ in
$L_{\La_1}\otimes\dots\otimes L_{\La_n}$, where $L_\La$ is the
irreducible $sl_{N+1}$ representation of highest weight $\La$.

In the present paper we propose a way to attack this
conjecture. Given an $h$-critical point with a given weight at infinity
$\La_\infty$, we describe a
procedure of constructing a set of $h$-critical points 
with weights at infinity of the form $w\cdot \La_\infty$, where $w$ is
in the $sl_{N+1}$ Weyl group and the dot denotes the shifted action,
see \Ref{shifted}. We call this procedure the reproduction procedure
and the resulting set of $h$-critical points a population.

The reproduction procedure makes use of a reformulation of the
algebraic system \Ref{Bethe slN} 
in terms of difference equations of the second order, and
furthermore in the following fertility property. Under some technical
conditions, zeroes of an $N$-tuple of 
polynomials $(y_1,\dots,y_N)$ form an $h$-critical
point if and only if there exist polynomials $\tilde
y_1,\dots,\tilde y_N$ such that 
\be
W(y_i(x),\tilde y_i(x))=y_{i-1}(x+h)y_{i+1}(x)T_i(x), 
\ee
where $W(u,v)=u(x+h)v(x)-u(x)v(x+h)$ is the discrete Wronskian and
polynomials $T_i$ are given explicitly in terms of $z_i$ and $\La_i$,
see \Ref{T polyn}. Furthermore, it turns out that in this case 
the zeros of the
tuple $(y_1,\dots,\tilde y_i,\dots,y_N)$ also form an $h$-critical
point and therefore we are able to repeat the same argument.
Thus, we get a family $h$-critical points, each is represented by
an $N$-tuple of polynomials.

The space $V$ of the first coordinates of these $N$-tuples has dimension
$N+1$ and it is called the fundamental space. 
The population is identified with the variety of all full flags in
$V$. Given a flag $\{0=F_0\subset F_1\subset\dots\subset F_{N+1}=V\}$, the
corresponding $N$-tuple of polynomials is given by $y_i=W_i(F_i)/U_i$,
where $U_i$ are some explicit polynomials written in terms of $T_i$,
see Lemma \ref{U=T},
and $W_i$ denotes the discrete Wronskian of order $i$. 

We consider the fundamental space of a population as an
$(N+1)$-dimensional subspace of the space $\C_d[x]$ of polynomials of
sufficiently big degree. Therefore a fundamental space defines a point in
the Grassmannian variety of $(N+1)$-dimensional subspaces of $C_d[x]$.
The fact that Wronskians of all $i$-dimensional subspaces in $V$ are
divisible by $U_i$ results in the conclusion that $V$ belongs to the
intersection of suitable Schubert cells in the Grassmannian. These
Schubert cells are related to special flags $\mc F(z_i)$ in
$\C_d[x]$. These flags are formed by the subspaces $F_j(z_i)$ 
in $\C_d[x]$ which consist of all polynomials divisible by
$(x-z_i)(x-z_i-h)\dots(x-z_i-(d-j)h)$. In addition,
the point of the Grassmannian $V$ belongs to a specific Schubert cell
related to the flag $\mc F(\infty)$. The flag $\mc F(\infty)$ is
formed by the subspaces $F_j(\infty)$ in $\C_d[x]$ which consist of 
all polynomials of degree less than $j$.

Any population contains at most 
one $h$-critical point associated to an integral dominant weight at infinity.
The fundamental spaces corresponding to different populations are
different. Therefore we related the problem of counting the
$h$-critical points associated to integral dominant weights at infinity
to the problem of counting the points in intersections of
Schubert cycles. It is well-known that
these Schubert cycles have the algebraic index of intersection equal to
the multiplicity of $L_{\La_\infty}$ in
$L_{\La_1}\otimes\dots\otimes L_{\La_n}$. It is an open question 
if for generic $z_i$ the intersection points are all of multiplicity one.

The scheme described above repeats the method in \cite{MV}, where a similar
picture was developed for the case of Bethe equations related to the
Goden model. The paper \cite{MV} is a smooth version of the present
paper; it uses Fuchsian differential equations and the usual Wronskians.

In the smooth case, the critical points, reproduction procedure and
populations are defined for any Kac-Moody algebra. It is not clear how
to achieve such level of generality
in the discrete situation of the present paper. We suggest
the relevant definitions in the case of root systems of type 
$B_N$ and $C_N$.

We say that 
zeroes of polynomials $(y_1,\dots,y_N)$ constitute a $B_N$ (resp. $C_N$)
$h$-critical point if zeroes of
\be
(y_1(x),\dots, y_N(x),y_{N-1}(x+h),y_{N-2}(x+2h),\dots,y_1(x+(N-1)h))
\ee
(resp.
\be
(y_1(x),\dots,y_N(x),y_N(x+h/2),y_{N-1}(x+3h/2),\dots,y_1(x+N-h/2)))
\ee 
form an $sl_{2N}$ (resp. $sl_{2N+1}$) $h$-critical point, see Section
\ref{def B section} (resp. Section \ref{def C section}) for details.
Our definition is motivated by similar properties of 
critical points in the smooth case,
see \cite{MV}. It turnes out that $h$-critical points defined in such
a way are exactly solutions of some systems of algebraic equations in both
cases of $B_N$
and $C_N$. We describe these algebraic equations.


A space of polynomials $V$ of dimension $N+1$ is called $h$-selfdual if the
space of discrete Wronskians 
$W(F)/U$, where $F$ runs over
$N$-dimensional subspaces of $V$ and $U$ is the greatest common
divisor of all $W(F)$, coincides with the space of functions of the form
$v(x-(N-1)h/2)$, where $v\in V$. Such a space $V$ has a natural
non-degenerate form which turns out to be symmetric if the dimension
of $V$ is odd and 
skew-symmetric if the dimension of $V$ is even.

A  population of $B_N$ (resp. $C_N$) type 
is then naturally identified with the 
space of isotropic flags in an $h$-selfdual space of dimension $2N$
(resp. $2N+1$). In particular a $B_N$ (resp. $C_N$) population is
identified with the $C_N$ (resp. $B_N$) flag variety.

The present paper deals with the $h$-analysis and additive
shifts. Similar results can be obtained in the case of the $q$-analysis
with multiplicative shifts related to the Bethe equations of the XXZ type.

The paper is constructed as follows. We start with the case of $sl_2$ in 
Section \ref{sl2}. The Sections \ref{slN} and \ref{slN fund} are
devoted to the case of $sl_{N+1}$. 
We discuss $h$-selfdual spaces of polynomials
in Section \ref{selfdual sec} and then deal with the cases of $B_N$ and $C_N$
in Sections \ref{B sec} and \ref{C sec}. 
Appendix A describes in detail the simplest
example of an $sl_3$ and ``$C_1$'' populations.  
Appendix B collects
identities involving discrete Wronskians.

\medskip 

We thank V. Tarasov for interesting discussions.

\section{The case of $sl_2$}\label{sl2}

\subsection{The Bethe equation}
Our parameters are distinct complex numbers $z_1,\dots,z_n$,
positive integers $\La_1,\dots,\La_n$ and a
complex non-zero number $h$. 
We call $z_i$ the {\it ramification points}, $\La_i$ the {\it weights} and $h$
the {\it step}.

For given $l\in\Z_{\geq 0}$, 
the $sl_2$ {\it Bethe equation} is the following system of algebraic
equations for the complex variables $\bs t=(t_1,\dots,t_l)$
\bean \label{Bethe sl2}
\prod_{s=1}^n \frac {t_j-z_s + \La_s h}
{t_j-z_s} \
\prod_{k \neq j}
\frac {t_j - t_k - h}
{t_j - t_k + h} \
 \ = \  1, 
\eean
where $j = 1, \dots , l$.
If all $\La_i$ are equal to $1$ then the system 
\Ref{Bethe sl2} is the Bethe equation
for the inhomogeneous XXX model, \cite{BIK}.

Consider the product of Euler gamma functions
\be
\Phi=\prod_{j=1}^l\prod_{s=1}^n \frac {\Gamma((t_j-z_s + \La_s h)/p)}
{\Gamma((t_j-z_s)/p)} \
\prod_{j,k; \;j<k}
\frac {\Gamma((t_j - t_k - h)/p)}
{\Gamma((t_j - t_k + h)/p)}.
\ee
The function $\Phi$ is called the master function associated to the
Yangian of $sl_2$  and $\bs z,\bs \La$. 
It is used in the construction of integral solutions to the
rational quantum Knizhnik-Zamolodchikov equation, see \cite{TV},
\cite{MV2}. 
The function
$\Phi$ is the
$sl_2$ difference counterpart of the master function (2.1) in \cite{MV}.
Equation \Ref{Bethe sl2} can be rewritten in the form
\be
\lim_{p\to 0}
\frac{\Phi(t_1,\dots,t_j+p,\dots,t_l)}{\Phi(t_1,\dots,t_j,\dots,t_l)}=1,\qquad
j=1,\dots,l. 
\ee

If $\bs t=(t_1,\dots,t_l)$ is a solution of \Ref{Bethe sl2} then any
permutation of the coordinates is also a solution.
An $S_l$ orbit of solutions of equation \Ref{Bethe sl2} such that
$t_i\neq t_k$ and
$t_i\neq z_s-jh$ for all $i,s,k$, $i\neq k$, and $j=1,\dots,\La_s$ is
called an 
{\it $h$-critical point}. 

Not all solutions of the Bethe equation are $h$-critical points.
For instance, if $n=0$ and $l$ is even, then  $t_1=\dots=t_l=0$ is a
solution of the Bethe equation which is not an $h$-critical point.

\subsection{Second order difference equations} 
For each $h$-critical point $\bs t$ we write the
monic polynomial $u(x)=\prod_{i=1}^l(x-t_i)$ and say that a polynomial
$cu(x)$ for any non-zero complex number $c$,
represents $\bs t$. 
In this section we 
rewrite the Bethe equation in terms of a difference equation for $u(x)$.

For two functions $u(x), v(x)$, define the discrete Wronskian
$W(u,v)(x)$ by the formula
\be
W(u,v)(x) = u(x+h)v(x) - u(x)v(x+h).
\ee
Let $A(x), B(x), C(x)$ be given functions. Consider the difference equation
\bean\label{eqn 1}
A(x)\ u(x+h) +B(x)u(x) +C(x) u (x-h) = 0 
\eean
with respect to the unknown function $u(x)$. Note that 
the complex vector space of polynomial solutions of equation \Ref{eqn 1} has
dimension at most 2.

The following lemma is straightforward.
\begin{lem}
If $u(x)$ and $v(x)$ are two solutions of equation \Ref{eqn 1}, then
the discrete Wronskian $W(u,v)(x)$ satisfies the first order
difference equation 
\bean\label{eqn wron}
W(u,v)(x+h) = \frac{C(x+h)}{A(x+h)} W(u,v)(x) .
\eean
\hfill $\square$
\end{lem}

Given ramification points $z_i$ and weights $\La_i$ we now fix the choice
of the functions
$A(x)$ and $C(x)$ as follows
\bean\label{AC}
A(x)=\prod_{s=1}^n (x - z_s), \qquad C(x)=\prod_{s=1}^n (x - z_s + \La_s h).
\eean
Then the polynomial solutions of 
equation \Ref{eqn wron} are constant multiples of the function 
\bean\label{eqn 3}
T(x) = \prod_{s=1}^n \prod_{i=1}^{\La_s}  (x - z_s + i h )  .
\eean


We say that a polynomial 
$u(x)$ is {\it generic} with respect to   $\bs z=(z_1,\dots,z_n)$,
$\bs \La=(\La_1,\dots,\La_n)$ if the polynomial $u(x)$ has no
multiple roots and no common roots with polynomials $u(x+h)$, $T(x)$.

\begin{lem}\label{bethe=diff}
Assume that $B(x)$ is a polynomial and 
that a solution $u(x)$ of \Ref{eqn 1} 
is a polynomial generic with respect to $\bs z, \bs \La$. 
Then $u(x)$ represents an $h$-critical point.

Conversely, if $u(x)$  represents an $h$-critical point then
the polynomial $u(x)$ satisfies equation \Ref{eqn 1}, where
$B(x)$ is a polynomial given by 
\bean\label{eqn 4}
B(x) =
- \left({\prod_{s=1}^n (x - z_s) 
u(x+h) + \prod_{s=1}^n (x - z_s + \La_s h) u (x-h) }\right)/
{u(x)} .
\eean
\end{lem}
\begin{proof}
The lemma follows from the fact that 
equation \Ref{Bethe sl2} with respect to $t_j$ is obtained by the
substituting $t_j$ in equation \Ref{eqn 1}. 
\end{proof}

\subsection{The $sl_2$ reproduction} 
By Lemma \ref{bethe=diff}, each $h$-critical point defines a
difference equation together with a polynomial solution.
This difference equation is a linear difference equation  
of the second order with polynomial coefficients.
In this section we show that 
the space of solutions of that difference equation contains
a complex two-dimensional space of polynomials
all of which (except for scalar multiples of finitely many) represent
$h$-critical points. 

\begin{lem}\label{two poly}
Let the roots of the polynomial $u(x)$ be distinct and 
satisfy the Bethe equation
\Ref{Bethe sl2}. Then the corresponding difference equation 
\Ref{eqn 1} has two linearly independent polynomial solutions.
\end{lem}
\begin{proof} 
Let $u(x)=\prod_{i=1}^n(x-t_i)$.
We look for the second polynomial solution $v(x)$ in the form  
$v(x)=c(x) u(x)$.
After multiplication by a
non-zero number, we have the equation
\bean\label{eqn 5}
v(x+h) u(x) -  v(x) u(x+h) = 
 \prod_{s=1}^n \prod_{i=1}^{\La_s}  (x - z_s + i h )  .
\eean
 Therefore $c(x)$ is a solution of the equation
\bean\label{eqn 6}
c(x+h)  - c(x) = \frac{ \prod_{s=1}^n \prod_{i=1}^{\La_s}  (x - z_s + i h )}
{u(x) u(x+h)} .
\eean
The right hand side of this equation is of the form
\be
f(x) + \sum_{j=1}^l \left(\frac{a_j}{x - t_j + h} - \frac{b_j}{x -
    t_j}\right),
\ee
where $f(x)$ is a polynomial and $a_j, b_j$ are some numbers.
Moreover, equation \Ref{Bethe sl2} imply that $a_j = b_j$ for all $j$.

The equation $c(x+h) - c(x) = f(x)$ has a polynomial solution. 
Denote it $\tilde c(x)$. Then $c(x) = \tilde c(x)
+ \sum_{j=1}^l {a_j}/({x - t_j})$ is a solution to \Ref{eqn 6}.

Then the function $v(x) = c(x) u(x)$ is a 
polynomial solution of equation \Ref{eqn 1}.
\end{proof}
Thus, starting from a solution $\bs t$ of the Bethe equation \Ref{Bethe
  sl2}, we obtain a difference equation of order two which has
a two-dimensional space of polynomial solutions.  
We call the projectivisation of this space the
 {\it  population of $h$-critical points related to $\bs t$} and denote
  $P(\bs t)$. 

Generic polynomials  with respect to $\bs z,\bs \La$ form a Zariski open
subset of the population. Roots of a generic polynomial form an $h$-critical
point by Lemma  \ref{bethe=diff}.

Let $\partial : f(x) \to f(x+h)$ be the shift operator acting on
functions of one variable.
The difference operator $\partial^2+B(x)/A(x)\partial+C(x)/A(x)$
obtained from an $h$-critical point 
$\bar{\bs t}\in P(\bs t)$ does not depend on the choice of $\bar{\bs
  t}$. We call this operator
the {\it fundamental operator associated to the population $P(\bs t)$} and
denote it by $D_{P(\bs t)}$.

\subsection{Fertile polynomials}
We give another equivalent condition for a polynomial to define an
$h$-critical point.

\begin{lem}\label{fertile}
Let a polynomial $u(x) =\prod_{i=1}^l (x - t_i)$ 
be generic with respect to $\bs z, \bs \La$.
Then the roots of $u(x)$ form a solution of equation \Ref{Bethe sl2}
if and only if 
there exists a polynomial solution $v(x)$ to equation \Ref{eqn 5}.
\end{lem}
\begin{proof} If the roots of $u(x)$ form a solution
of equation  \Ref{eqn 3}, then equation \Ref{eqn 1} with the polynomial
$B(x)$ given by
\Ref{eqn 4} has two linearly 
independent polynomial solutions by Lemma \ref{two poly}, and hence
equation \Ref{eqn 5} has a polynomial solution. 

Now assume that $v(x)$ is
a polynomial solution to \Ref{eqn 5}. Then the rational
function $c(x) = v(x)/u(x)$ is a solution of equation \Ref{eqn 6}. Then
\be
\text{Res}_{x=t_j}  \frac{ \prod_{s=1}^n \prod_{i=1}^{\La_s}  (x - z_s + i h )}
{u(x)\ u(x+h)}  + 
\text{Res}_{x=t_j-h}  \frac{ \prod_{s=1}^n \prod_{i=1}^{\La_s} (x - z_s + i h )}
{u(x) u(x+h)}= 0
\ee
for $j = 1, \dots , l$. This system of equations is equivalent to equation
\Ref{Bethe sl2}.
\end{proof}

A polynomial $u(x)$ is called {\it fertile} if there exist a polynomial
$v(x)$ such that the discrete Wronskian of $u$ and $v$ is $T(x)$, 
$W(u,v)(x)=T(x)$. All polynomials
in a population of $h$-critical points are fertile.
Moreover, the $h$-critical points correspond to fertile generic polynomials. 

We get the following immediate corollary.
\begin{lem}\label{elem W}
If equation \Ref{eqn 1} has two polynomial solutions
of degree $l$ and $l'$, $l \neq l'$,
then $l + l' - 1 = \La_1 + \dots + \La_n$.
\hfill $\square$
\end{lem}

\section{The case of $sl_{N+1}$}\label{slN}

Let $\al_1,\dots,\al_N$ be the simple roots of $sl_{N+1}$. We have
$(\al_i,\al_i)=2$, $(\al_i,\al_{i\pm 1})=-1$ with all other scalar
  products equal to zero. 

\subsection{Definition of $h$-critical points}
In the case of $sl_{N+1}$ we fix the following parameters:
{\it ramification points} $\bs z=(z_1,\dots,z_n)\in\C^n$; 
non-zero dominant integral $sl_{N+1}$
{\it weights} $\bs \La=(\La_1,\dots,\La_n)$, {\it relative shifts}
 $\bs b^{(i)}=(b^{(i)}_1,\dots,b^{(i)}_n)\in\C^n$, $i=1,\dots,N$, and the
 step $h\in\C$, $h\neq 0$. 

We call the set of parameters $\bs z, \bs \La, \bs b^{(i)}$,
$i=1,\dots,N$, the {\it 
initial data}. We denote $\La_s^{(i)}=(\La_s,\al_i)\in\Z_{\geq 0}$.

Let $\bs l = (l_1, \dots , l_N) \in \Z^N_{\geq 0}$.
The $sl_{N+1}$-weight 
\bean\label{w at inf}
\La_\infty  =  \sum_{s=1}^n \La_s  - \sum _{i=1}^N  l_i  \al_i,
\eean
is called the {\it weight at infinity}.

The following system of algebraic equations for variables 
$\bs t = ( t^{(i)}_j)$, $i=1,\dots,N$, $j = 1,\dots ,l_i$,
\begin{align}
\prod_{s=1}^n \frac {t_j^{(i)}-z_s + b_s^{(i)} h + \La_s^{(i)} h}
{t_j^{(i)}-z_s + b_s^{(i)} h} \
\prod_{k=1}^{l_{i-1}}
\frac{t_j^{(i)} - t_k^{(i-1)} + h}
{t_j^{(i)} - t_k^{(i-1)}} \
& \prod_{k \neq j}
\frac {t_j^{(i)} - t_k^{(i)} - h}
{t_j^{(i)} - t_k^{(i)} + h} \times \notag \\
& \times \prod_{k=1}^{l_{i+1}}
\frac{t_j^{(i)} - t_k^{(i+1)}}
{t_j^{(i)} - t_k^{(i+1)} - h} =  1, \label{Bethe slN}
\end{align}
is called the
$sl_{N+1}$ {\it Bethe equation} associated with the initial data and
the weight at infinity.

Note that by a shift of variables $\tilde t_i^{(j)}= t_i^{(j)}-jh/2$
the system can be written in a slightly more symmetrical way.

Note that in the quasiclassical limit $h\to 0$, system \Ref{Bethe
  slN} becomes system (2.2) of \cite{MV} specialized to the case of 
$sl_{N+1}$. 


The product of symmetric groups
$S_{\bs l}=S_{l_1}\times \dots \times S_{l_N}$ acts on the  set of
solutions of \Ref{Bethe slN}
permuting the coordinates with the same upper index. 
An $S_{\bs l}$ orbit of solutions of the Bethe equation such that
$t_j^{(i)}\neq t_k^{(i)}$,
$t_j^{(i)}\neq z_s-b_s^{(i)}-rh$ for 
 all $j,i,s,k$, $k\neq j$ and $r=1,\dots,\La_s^{(i)}$
 is called an {\it $h$-critical point associated with the initial data
   and the weight at infinity}.

Let $L_\La$ be the irreducible $sl_{N+1}$ module with highest weight $\La$.
Let the initial data satisfies $b_s^{(j)}=-\sum_{i=1}^j\La_s^{(i)}$.
We have the following conjecture.

\begin{conj} 
If $\La_\infty$ is integral dominant then for generic $\bs z$ the
number of $h$-critical points associated with the initial data and the 
weight at infinity equals to the multiplicity of $L_{\La_{\infty}}$ in
$L_{\La_1}\otimes\dots \otimes L_{\La_n}$.
\end{conj}

\subsection{Difference equations of the second order}
In this section, given an $h$-critical point we obtain $N$ 
difference operators of
order two. The $i$th operator is the fundamental operator with respect to
the $sl_2\in sl_{N+1}$ in the direction $\al_i$.

Let $\bs t$ be an $h$-critical point.
For $i=1,\dots, N$, introduce polynomials 
\bean
y_i(x)&=&\prod_{j=1}^{l_i}(x - t^{(i)}_j).\label{y polyn}
\eean
We also set $y_0(x) = 1$,  $y_{N+1}(x) = 1$.

The $N$-tuple $\bs y$ uniquely determines the $h$-critical point $\bs t$
and we say that  $\bs y$ {\it represents $\bs t$}. 
We consider the tuple $\bs y$ up to multiplication 
of each coordinate by a non-zero number, 
since we are interested only in the roots of polynomials $y_1, \dots ,
y_N$. Thus the tuple defines a point in the direct product
$\PCN$ of $N$ copies of the projective space associated with the vector 
space of polynomials of $x$. In this paper we view all $N$-tuples
of polynomials as elements of $\PCN$.

Introduce the polynomials
\bean 
T_i(x)&=& \prod_{s=1}^n \prod_{j=1}^{\La_s^{(i)}} (x - z_s + b_s^{(i)} h + j
h),\label{T polyn} \\
A_i(x) &=& 
\prod_{s=1}^n (x - z_s + b_s^{(i)} h)
\prod_{k=1}^{l_{i-1}} (x - t_k^{(i-1)})
\prod_{k=1}^{l_{i+1}} (x - t_k^{(i+1)} - h),\notag \\
C_i(x) &=& 
\prod_{s=1}^n (x - z_s + b_s^{(i)} h + \La_s^{(i)} h) 
\prod_{k=1}^{l_{i-1}} (x - t_k^{(i-1)} + h) 
\prod_{k=1}^{l_{i+1}} (x - t_k^{(i+1)}) .\notag
\eean

We say that a tuple of polynomials $\bs y \in \PCN$ is {\it generic
  with respect to the initial data} if for all $i$ the polynomial
  $y_i(x)$ has no multiple roots and no common roots with polynomials
$y_i(x+h)$, $y_{i-1}(x+h)$, $y_{i+1}(x)$, $T_i(x)$.
 
For $i=1,\dots,N$, consider difference equations of the second
order of the form
\bean\label{eqn i}
A_i(x)  u ( x + h )
+ B_i(x) u ( x )  +  C_i(x) 
 u ( x - h ) =0,
\eean
where $B_i(x)$ are any functions.

\begin{lem}\label{bethe=diff slN}
Assume that for all $i$, $y_i$ is a 
polynomial solution of $i$-th equation in \Ref{eqn i} for some
polynomial $B_i(x)$. If the tuple $\bs
y=(y_1,\dots,y_N)$ is generic with respect to the initial data then 
$\bs y$ represents an $h$-critical point.

Conversely, if a tuple of polynomials $\bs y=(y_1,\dots,y_N)$
represents an $h$-critical point, then for  all $i$ the polynomial
$y_i$ satisfies the $i$-th equation in \Ref{eqn i}, where
$B_i(x)$ is a polynomial given by 
$B_i(x) =-(A_i(x)y_i(x+h)+C_i(x)y_i(x-h))/y_i(x)$.
\end{lem}
\begin{proof}
The lemma follows directly from the $sl_2$ counterpart,
Lemma \ref{bethe=diff}. \end{proof}

\subsection{Fertile tuples}
In this section we 
discuss yet another criteria for a tuple of polynomials to represent an
$h$-critical point.

We say that a tuple of polynomials $\bs y \in \PCN$ is {\it fertile
  with respect to the initial data} 
if for every $i$ there exists a polynomial $\tilde y_i$
such that the discrete Wronskian is 
\bean\label{main}
W(y_i , \tilde y_i)(x)  = T_i(x)   y_{i-1}(x+h) y_{i+1} (x) .
\eean

\begin{lem}\label{fertile=critical}
A generic tuple $\bs y$ is fertile if and only if it represents
an $h$-critical point. 
\end{lem}
\begin{proof}
The lemma follows from the $sl_2$ considerations, 
see Lemma \ref{fertile}. \end{proof}

Let $\bs y$ be fertile and let $y_i, \tilde y_i$ satisfy \Ref{main}. Then
polynomials of the form $c_1 y_i + c_2 \tilde y_i$ with $c_1, c_2 \in \C$
 span the two dimensional space of polynomial solutions of equation
 \Ref{eqn i}. 
The tuples $(y_1,\dots,c_1 y_i + c_2 \tilde y_i,\dots,y_N)\in\bs
P(\C[x])^N$ are called  
{\it immediate descendents} of $\bs y$ in the $i$-th direction.
Note that if $\bs y$ is generic, then all but finitely many immediate
descendents are generic.

Let $\C_d[x]$ be the space of polynomials of degree not greater than $d$.
The set of fertile tuples is closed in $\bs P (\C_d[x])^N$:

\begin{lem}\label{limit fertile lem}
Assume that a sequence  of fertile tuples of 
polynomials $\bs y_k$, $k = 1, 2, \dots$, 
has a limit $\bs y_\infty$ in $\bs P (\C_d[x])^N$ as $k$
tends to infinity. Then the limiting tuple $\bs y_\infty$ is fertile.
Moreover, if $\bs y^{(i)}_\infty$ is 
an immediate descendant of $\bs y_\infty$ 
in the $i$-th direction, then
there exist immediate descendants
$\bs y_k^{(i)}$ of $\bs y_k$ in the $i$-th direction such that 
$\bs y^{(i)}_\infty$ is the limit of $\bs y_k^{(i)}$.
\end{lem}
\begin{proof}
Let us prove $\bs y_\infty$ is fertile in direction $i$. For each $k$
we have a difference equation \Ref{eqn i} which we denote $E_k$ 
and a plane of polynomial
solutions $P_k$. By the assumptions of the lemma there is a limiting
equation $E_\infty$ as $k\to \infty$.

The space $P_k$ defines a point in the projective
Grassmanian variety of planes in  $\C_d[x]^N$. 
Let $P_\infty$ be a limiting plane, then all polynomials in $P_\infty$ 
are solutions of equation $E_\infty$  and the lemma is proved. 
\end{proof}

\subsection{Reproduction procedure}
Given an $h$-critical point, we describe a procedure of obtaining a
family of $h$-critical points.

\begin{theorem}
Let $\bs y$ represent an $h$-critical point and let 
$y_i, \tilde y_i$ satisfy \Ref{main}.
Assume that $\bs y^{(i)}=(y_1,\dots,\tilde y_i,\dots,y_N)$ is generic. Then $\bs y^{(i)}$ represents an 
$h$-critical point.
\end{theorem}
\begin{proof} 
Denote $\bar t^{(i)}_j$ the roots of $\bar y_i = c_1 y_i + c_2 \tilde y_i$. 
Let $z$ be a root of $T_i(x) y_{i-1}(x+h) y_{i+1}(x)$. Then we have
\be
\frac{y_i(z+h)}{y_i(z)} = \frac{\bar y_i(z+h)}{\bar 
y_i(z)}.
\ee
In particular we choose $z = t^{(i-1)}_k-h$ and $z = t^{(i+1)}_k$ and obtain
\bea
\prod_j \frac { t^{(i-1)}_k -  t^{(i)}_j}{ t^{(i-1)}_k -  t^{(i)}_j - h}
&=&
\prod_j \frac { t^{(i-1)}_k -  \bar t^{(i)}_j}{ t^{(i-1)}_k -  \bar
  t^{(i)}_j - h} , 
\\
\prod_j \frac { t^{(i+1)}_k -  t^{(i)}_j + h}{ t^{(i+1)}_k -  t^{(i)}_j}
&=&
\prod_j \frac { t^{(i+1)}_k -  \bar t^{(i)}_j + h}{ t^{(i+1)}_k -
  \bar t^{(i)}_j}  
\eea
for all $k$ and $i$. In addition Lemma \ref{bethe=diff} implies that
for all $k$,
\be
\prod_{j;\ j\neq k}\frac{t^{(i)}_k -  t^{(i)}_j - h}{ t^{(i)}_k -
    t^{(i)}_j+h}=\prod_{j;\ j\neq k}\frac{ \bar t^{(i)}_k - \bar
      t^{(i)}_j - h}{\bar t^{(i)}_k - \bar t^{(i)}_j+h}.
\ee
That proves the theorem.
\end{proof}
Thus, starting with a tuple $\bs y$,
representing an $h$-critical point and an index 
$i \in \{1, \dots , N\}$, we  construct a closed subvariety of immediate
descendants $\bs y^{(i)}$ in $\PCN$ which is isomorphic to $\bs P^1$.
All but finitely many points in this subvariety correspond to
$N$-tuples of polynomials which 
represent $h$-critical points. All constructed $N$-tuples of
polynomials are fertile. We call this construction the {\it
  simple reproduction procedure in the $i$-th direction}. 

Now we can repeat the simple reproduction procedure in some other 
$j$-th direction applied to all points obtained in the previous step
and obtain more critical points. 

We continue the process until no simple reproduction in any direction 
applied to any polynomial obtained in the previous steps produces 
a new polynomial. The result is called the {\it population originated at the
$h$-critical point $\bs y$} and is denoted $P(\bs y)$.

More formally, the population of $\bs y$ is the set of all $N$-tuples
of polynomials $\bar {\bs y}\in\PCN$ such that there exist $N$-tuples $\bs
y_1=\bs y$, $\bs y_2,\dots , \bs y_k=\bar{\bs y}$, such that $\bs y_i$
is an immediate descendant of $\bs y_{i-1}$ for all $i=2,\dots,k$.

Obviously, if two populations intersect then they coincide.  


\section{Fundamental space of an $sl_{N+1}$ population}\label{slN fund}
Given an $sl_{N+1}$ population of critical points we construct a space of
polynomials called the fundamental space. 
Elements of the population are in natural correspondence
with full flags in the fundamental space. 
Then we show that the problem of counting $h$-critical points is
equivalent to a problem of Schubert calculus.

\subsection{Degrees of polynomials in a population and the $sl_{N+1}$ Weyl group}
Let $\h$ be the Cartan subalgebra of $sl_{N+1}$ and let 
$(\ , \ )$ be the standard scalar product on $\h^*$.
Denote $\rho \in \h^*$ the half sum of positive roots.

The Weyl group $\mathcal W \in \text{End} (\h^*)$ is generated by 
reflections $s_i$, $i = 1, \dots , N$, 
\be
s_i(\la) = \la - (\la,\al_i)  \al_i , \qquad  \la \in \h^* .
\ee
We use the notation
\bean\label{shifted}
w \cdot \la = w( \la + \rho) -  \rho,
\qquad w \in \mathcal W,\; \la \in \h^* ,
\eean
for the shifted action of the Weyl group.

Let an initial data be given.
Let $\bs y = (y_1, \dots , y_N) \in \PCN$, and
let $l_i $ be the degree of the polynomial $y_i$ for $i = 1, \dots N$.
Let $\La_\infty$ be the weight at infinity defined by \Ref{w at inf}.

Assume the tuple $\bs y$ is fertile.
Let the tuple $ {\bs y}^{(i)} = (y_1, \dots , \tilde y_i, \dots , y_r)$ be 
an immediate descendant of $\bs y$ in the direction $i$ and
let  $ \La^{(i)}_\infty$ be the $sl_{N+1}$ weight at infinity of 
${\bs y}^{(i)}$.

The following lemma follows from Lemma \ref{elem W}.

\begin{lem}\label{simple inf lem} 
If the degree of $\tilde y_i$ 
is not equal to the degree of  $y_i$, then  
\be
\Lambda^{(i)}_\infty\ = \ s_i \cdot {\La}_\infty\ ,
\ee
where $s_i\cdot $ is the shifted action of the $i$-th generating
reflection of the Weyl group.
\end{lem}

Therefore we obtain the following 

\begin{prop}\label{inf weight thm}
Let a tuple $\bs y^0$ represent an $h$-critical point associated with
a given initial data and weight at infinity $\La_\infty$.
Let $P$ be the population of $h$-critical points originated at  $\bs y^0$.
Then
\begin{enumerate}
\item[$\bullet$]
For any tuple $\bs y \in P$,  there is an element $w$ of the
$sl_{N+1}$ Weyl group
$\mathcal W$, 
such that the weight at infinity of $\bs y$
is  $w \cdot \Lambda_\infty$.
\item[$\bullet$] 
For any element $w \in  \mathcal W$,  there is a tuple $\bs y \in P$ 
whose weight at infinity
is  $w \cdot \Lambda_\infty$.
\end{enumerate}
\hfill $\square$
\end{prop}

Proposition \ref{inf weight thm} gives sufficient conditions for absence
of $h$-critical points.

\begin{cor}
There is no $h$-critical point
in either of the two cases 
\begin{enumerate}
\item
if there is an element $w$ 
of the Weyl group such that $\sum_{s=1}^n \La_s - w \cdot  \La_\infty$ 
does not belong to the cone $\Z_{\geq 0} \al_1 \oplus \dots \oplus  
\Z_{\geq 0} \al_N$;
\item if $\La_\infty$ belongs to one of the reflection hyperplanes 
of the shifted action of the Weyl group.
\end{enumerate}
\end{cor}

The next corollary says that under certain conditions on weights there is 
exactly one population of $h$-critical points.

The tuple $(1, \dots , 1) \in \PCN $ is the unique $N$-tuple of
non-zero polynomials of degree 0. 
The weight at infinity of $(1, \dots , 1)$ is
$\Lambda_{\infty, (1, \dots , 1)} 
= \sum_{s=1}^n \La_s$. 
Let $P_{(1, \dots , 1)}$ be the population associated  to the initial data 
and originated at $(1, \dots , 1)$.

\begin{cor}
Let $\bs y$ be an $h$-critical point such that
$\La_\infty$ has the form $w\cdot\Lambda_{\infty, (1, \dots , 1)}$ for some
$w\in\mathcal W$. Then $\bs y$ belongs to the population  
$P_{(1, \dots , 1)}$.
\end{cor}

\subsection{The difference operator of a  population}
In this section we describe a linear difference operator of order $N+1$
related to a population of $h$-critical points.

Let an initial data  $\bs z, \bs \La, \bs b^{(i)}$ be given.
Let the polynomials $T_i$, $i = 1, \dots , N$, be defined by formula 
\Ref{T polyn}.

 Let $\bs y =(y_1, \dots, y_N)$ be an $N$-tuple of non-zero polynomials. Set
 $y_0 = y_{N+1} = 1$. Define a  linear difference operator of order $N+1$ 
with meromorphic coefficients 
\bean\label{Diff}
&&D (\bs y)  =  (   \partial  - 
\frac{y_{N}(x)}{y_{N}(x+h)}
\prod_{s=1}^N \frac{T_{s}(x + (N-s+1)h)}{T_{s}(x + (N-s)h)} )\times
\notag
\\
&&
(   \partial  - 
\frac{y_{N}(x+h)}{y_{N}(x)} 
\frac{y_{N-1}(x)}{y_{N-1}(x+h)}
\prod_{s=1}^{N-1} \frac{T_{s}(x + (N-s)h )}{T_{s}(x + (N-1-s) h )}
 ) \times \dots\times
\notag
\\                 
& {} &
(   \partial  - 
\frac{y_{2}(x+h)}{y_{2}(x)} 
\frac{y_{1}(x)}{y_{1}(x+h)}
 \frac{T_{1}(x+h)}{T_{1}(x)} ) \ 
(   \partial  - 
\frac{y_{1}(x+h)}{y_{1}(x)} )=
\notag
\\
 &=&
\prod^{0 \to N}_i
\left( \partial  -  
\frac{y_{N+1-i}(x+h)}{y_{N+1-i}(x)} 
\frac{y_{N-i}(x)}{y_{N-i}(x+h)}
\prod_{s=1}^{N-i} \frac{T_{s}(x + (N-i-s+1)h)}{T_{s}(x+(N-i-s)h)}
 \right).
\eean
Notice that the first coordinate $y_1$ of the $N$-tuple  belongs to 
the kernel of the operator $D(\bs y)$.

\begin{theorem}\label{ind thm}
Let $P$ be a population of $h$-critical points originated at some tuple
$\bs y^0$. Then the operator $D(\bs y)$ does not depend on the choice
of $\bs y$ in $P$.
\end{theorem}
\begin{proof} We have to prove that if $\bs y, \tilde{\bs y}\in P$ 
then $D ( \bs y ) = D ( \tilde {\bs y} )$. It is enough to show that
is the case when $\tilde{\bs y}$ is an immediate descendant of $\bs y$
in some direction $i$ and both $\bs y, \tilde{\bs y}$ 
represent critical points. Therefore, we
assume  that $y_j = \tilde y_j$, for all $j,  j \neq  i$,  and
$W(y_i, \tilde y_i)(x)  =  T_i(x)\  y_{i-1}(x+h) y_{i+1}(x)$. 
In this case, all  factors of $D(\bs y)$ and
$D(\tilde{\bs y})$ are the same except the two  factors 
which involve $y_i$ or $\tilde y_i$. So it is left to show 
that for any function $u(x)$ we have
\bea
(   \partial & -&  
\frac{y_{i+1}(x+h)}{y_{i+1}(x)} 
\frac{y_{i}(x)}{y_{i}(x+h)}
\prod_{s=1}^{i} \frac{T_{s}(x+(i-s+1)h)}{T_{s}(x+(i-s)h)} )\times
\\
&&
(   \partial  -  
\frac{y_{i}(x+h)}{y_{i}(x)} 
\frac{y_{i-1}(x)}{y_{i-1}(x+h)}
\prod_{s=1}^{i-1} \frac{T_{s}(x+(i-s)h)}{T_{s}(x+(i-s-1)h)} )\ {} \ u(x) \ =
\\
(   \partial  &- & 
\frac{y_{i+1}(x+h)}{y_{i+1}(x)} 
\frac{\tilde y_{i}(x)}{\tilde y_{i}(x+h)}
\prod_{s=1}^{i} \frac{T_{s}(x+(i-s+1)h)}{T_{s}(x+(i-s)h)} )\times
\\
&&
(   \partial  -  
\frac{\tilde y_{i}(x+h)}{\tilde y_{i}(x)} 
\frac{y_{i-1}(x)}{y_{i-1}(x+h)}
\prod_{s=1}^{i-1} \frac{T_{s}(x+(i-s)h)}{T_{s}(x+(i-s-1)h)} )
\ {} \ u(x) \ .
\eea
After the change of variables 
\be
v(x)
 = u(x) \frac{\prod_{s=1}^{i-1}T_s(x+(i-s-1)h)}{ y_{i-1}(x)},
\ee
we have to prove that
 \bea
(   \partial  -  
\frac{y_{i+1}(x+h)}{y_{i+1}(x)} 
\frac{y_{i}(x)}{y_{i}(x+h)}
\frac{y_{i-1}(x+2h)}{y_{i-1}(x+h)} 
\frac{T_{i}(x+h)}{T_{i}(x)} )
(   \partial  -  
\frac{y_{i}(x+h)}{y_{i}(x)} )\  v(x)
\ =
\\
(   \partial  -  
\frac{y_{i+1}(x+h)}{y_{i+1}(x)} 
\frac{\tilde y_{i}(x)}{\tilde y_{i}(x+h)}
\frac{y_{i-1}(x+2h)}{y_{i-1}(x+h)} 
\frac{T_{i}(x+h)}{T_{i}(x)} )
(   \partial  -  
\frac{\tilde y_{i}(x+h)}{\tilde y_{i}(x)} )\  v(x) \ .
\eea
This identity is easily checked directly using the equation connecting $y_i$ and $\tilde y_i$.
\end{proof}
The difference operator \Ref{Diff} is called the {\it fundamental
  operator associated to the
  population $P$} and is denoted $D_P$.

\begin{cor}\label{first coor cor}
  Let $\bs y$ be a member of a population $P$.  Then the first
  coordinate $y_1$ of ${\bs y}$ is in the kernel of the operator
  $D_P$.
\end{cor}

\begin{prop}\label{DP prop}
The space of polynomial solutions to the difference equation 
\bean\label{DP eq}
D_P u = 0
\eean
has dimension $N+1$.
\end{prop}

\begin{proof}
Let $\bs y$ be a member of the population $P$. Assume that $\bs y$
represents an $h$-critical point. 
We construct polynomials $u_1, \dots , u_{N+1}$,
satisfying equation \Ref{DP eq}. 

Set $u_1 = y_1$. The polynomial $u_1$ is a solution of \Ref{DP eq}.

Let $u_2$ be a polynomial such that $W(u_1,u_2)(x) \ = \ T_1(x) y_2(x)$. 
The polynomial $u_2$ is a solution of \Ref{DP eq}.

Let $\tilde y_2$ be a polynomial such that 
$W ( y_2, \tilde y_2 )(x) \ = \ T_2(x) y_1(x+h) y_3(x)$. 
Such a polynomial  can be chosen so that
$( y_1, \tilde y_2 , \dots , y_N )$ is generic and therefore
represents an $h$-critical point.
Choose a polynomial $u_3$ to satisfy  equation
$W(u_3,y_1)(x) \ =\ T_1(x) \tilde y_2(x) $. 
The polynomial $u_3$ is a solution of \Ref{DP eq}.

In general, to construct a polynomial
$u_{i+1}$ we find $\tilde {\bs  y}_i = (y_1, \dots , \tilde y_{i}, \dots , y_N)$, such that 
$\tilde {\bs y}_i$ is
  generic  and $W(y_i, \tilde y_i)(x) = T_i(x) y_{i-1} (x+h) y_{i+1}(x)$ and then repeat
  the construction
 for $u_i$ using $\tilde{\bs y}_i$ instead of $\bs y$.

Let $V$ be the complex vector space spanned by polynomials $u_1, \dots ,
u_{N+1}$.  We show  that the space $V$ has dimension $N+1$.

Let $W(g_1,\dots,g_s$ be the discrete Wronskian of functions $g_1,
\dots, g_s$ defined by
\be\label{wronskian}
W (g_1, \dots , g_s)(x) = \det(g_i(x + (j - 1)h))_{i,j=1}^s
\ee
For $s=0$, we define the corresponding discrete
Wronskian to be $1$. 

\begin{lem}\label{wr(u) lem}
For $i = 1, \dots , N+1$, we have 
\be
W(u_1, \dots, u_i)(x)\ =\ y_i(x) \ \prod _{j=1}^{i-1}T_1(x+(j-1)h)\ 
\prod_{j=1}^{i-2} T_2(x+(j-1)h)\  \dots\  T_{i-1}(x) .
\ee
\end{lem}
\begin{proof}
The lemma is proved in the same way as Lemma 5.5 in \cite{MV} making
use of Wronskian identities in Lemmas \ref{f wronskian} and \ref{wr id 1}.
\end{proof}
The proposition is proved.
\end{proof}

A linear difference equation  of order $N+1$ 
cannot have more than $N+1$ polynomial solutions
linearly independent over $\C$. 
The complex $(N+1)$-dimensional vector space of
polynomial solutions of 
equation \Ref{DP eq} is 
called {\it the fundamental space of the population} $P$
and is denoted $V_P$.

\subsection{Frames} We describe an $h$-analogue of ramification
 properties of a space of polynomials.

A space of polynomials $V$ is called a {\it space without base points} if for
any $z\in\C$ there exists $v\in V$ such that $v(z)\neq 0$.

Let $V$ be an $(N+1)$-dimensional vector space of polynomials without
base points. 
Let $U_i(x)$ be the monic polynomial which is the greatest common  
divisor of the family of polynomials
$\{W(u_1,\dots,u_i)\ |\ u_1,\dots,u_i\in V\}$.

\begin{lem}\label{U=T}
There exist a unique sequence of monic polynomials
$T_1(x),\dots,T_N(x)$ such that 
\be
U_i(x)=\prod _{j=1}^{i-1}T_1(x+(j-1)h)\ \prod_{j=1}^{i-2}
T_2(x+(j-1)h)\  \dots\  T_{i-1}(x);
\ee
for $i=1,\dots,N+1$.
\end{lem}

We call the sequence of monic 
polynomials $T_1(x), \dots , T_N(x)$ the {\it frame} of $V$.

\begin{proof}
We construct the polynomials $T_i$ by induction on $i$. 
For $i=0$, we have $U_1=1$. 
For $i=1$ we just set $T_1=U_2$. Suppose the lemma is proved for all
$i=1,\dots,i_0-1$. Then we set
\be
S(x)=\prod_{i=1}^{i_0-2}\prod _{j=1}^{i_0-i}T_i(x+(j-1)h), \qquad 
T_{i_0-1}(x)=U_{i_0}(x)/S(x).
\ee
We only have to show that $T_{i_0-1}$ is a polynomial. 
In other words, we have to show that a Wronskian of any 
$i_0$ dimensional subspace in
$V$ is divisible by $S(x)$. 

Consider the Grassmanian $Gr(i_0-2,V)$ 
of $(i_0-2)$-dimensional spaces in $V$. For any $z\in \C$ the set of
points in $Gr(i_0-2,V)$ such that the corresponding Wronskian divided
by $U_{i_0-2}$ does not vanish at $z$, is a Zariski open algebraic
set. Therefore we have a Zariski open set of points in
$Gr(i_0-2,V)$ such that the
corresponding  Wronskian divided
by $U_{i_0-2}$ does not vanish at roots of $S(x-h)$. We call such
subspaces acceptable.

Therefore we have a Zariski open set of points in $Gr(i_0,V)$ such
that the corresponding $i_0$ dimensional space contains an acceptable 
$i_0-2$ dimensional subspace. Let $u_1,\dots,u_{i_0}\in V$ are such 
that $u_1,\dots,u_{i_0-2}$ span an acceptable space.
It is enough to show that $W(u_1,\dots,u_{i_0})$ is divisible by 
$S(x)$.

Then using Wronskian 
identities in Lemmas \ref{f wronskian} and \ref{wr id 1} we have for
suitable polynomials $f_1,f_2,g$:
\bea
W(u_1,\dots,u_{i_0})=\frac{W(W(u_1,\dots,u_{i_0-1}),W(u_1,\dots,u_{i_0-2},u_{i_0}))}{W(u_1,\dots,u_{i_0-2})(x+h)}=\\
\frac{W(U_{i_0-1}f_1,U_{i_0-1}f_2)}{U_{i_0-2}(x+h)g(x+h)}=
\frac{U_{i_0-1}(x)U_{i_0-1}(x+h)}{U_{i_0-2}(x+h)}\;\frac{W(f_1,f_2)}{g(x+h)}=S(x) \frac{W(f_1,f_2)}{g(x+h)}.
\eea
Since the space spanned by $u_1,\dots,u_{i_0}\in V$ is accepatble, 
the polynomial $g(x+h)=W(u_1,\dots,u_{i_0-2})(x+h)/U_{i_0-2}(x+h)$ and
$S(x)$ are relatively
prime. Therefore the Wronskian $W(u_1,\dots,u_{i_0})$ is divisible by 
$S(x)$.
\end{proof}

\subsection{Frames of the fundamental spaces}
Let $\bs y=(y_1,\dots,y_N)$ represent an $h$-critical point associated to
the initial data
$\bs z,\bs\La,\bs b^{(i)}$. Let $P$ be a population of $h$-critical
points originated at $\bs y$.

\begin{prop}\label{T form frame} The fundamental space $V_P$ of the
  population $P$ has no base points. 
The polynomials $T_1(x), \dots, T_N(x)$ given by \Ref{T polyn}
form a frame of $V_P$.
\end{prop}
\begin{proof}
We construct polynomials $u_1,\dots,u_{N+1}\in V_P$ as in the proof of
 Proposition \ref{DP prop}. We have
\bea
&&W(u_1,\dots,u_i)=y_i(x) \prod_{j=1}^{i-1}\prod_{r=1}^{i-j}T_j(x+(r-1)h), \\
&&W(u_1,\dots,u_{i-1},u_{i+1})=
\tilde y_i(x)\prod_{j=1}^{i-1}\prod_{r=1}^{i-j}T_j(x+(r-1)h),
\eea
where 
$(y_1,\dots,\tilde y_i,\dots, y_N)$ is a descendant of $\bs y$ in the $i$
direction. 

In particular $W(u_1,\dots,u_i)$ is divisible by 
$\prod_{j=1}^{i-1}\prod_{r=1}^{i-j}T_j(x+(r-1)h)$. Therefore
Wronskians of all $i$ dimensional planes are divisible by this polynomial.

Moreover, we have $W(y_i,\tilde y_i)=T_i(x)y_{i-1}(x+h)y_{i+1}(x)$. Since
  $y_i$ and $T_i(x)y_{i-1}(x+h)y_{i+1}(x)$ have no common roots, the
  polynomials $y_i(x)$ and $\tilde y_i(x)$ have no common roots. It
  follows that the greatest common divisor of $i$-dimensional
  Wronskians is $\prod_{j=1}^{i-1}\prod_{r=1}^{i-j}T_j(x+(r-1)h)$.

The absence of base points for $V_P$ follows from the case of $i=1$.
\end{proof}

The converse statement is also true.

\begin{prop} \label{space=popul} 
Let $V$ be a space of polynomials of dimension $N+1$
  without base points and with the frame $T_i$. Let $\bs z,\bs \La,\bs
  b^{(i)}$ be the initial data related to
  polynomials $T_i$ by \Ref{T
  polyn}. Then $V$ is the fundamental space of a population of
  $h$-critical points associated to the initial data $\bs z,\bs \La,\bs
  b^{(i)}$.
\end{prop}
\begin{proof}
We postpone the proof until after Section \ref{flag section}, where
the generating morphism $\beta$ is described. Then
the proof is similar to the proof of Proposition 5.17 in \cite{MV} 
combined with Lemma 5.20 in \cite{MV}.
\end{proof}

From Propositions \ref{T form frame} and \ref{space=popul} 
we obtain the following theorem.

\begin{theorem}
There is a bijective correspondence between populations of $h$-critical
points associated to a given initial data 
and the spaces of polynomials with framing $T_i$, where $T_i$ are
related to the initial data by \Ref{T polyn}. $\square$ 
\end{theorem}

\subsection{Schubert calculus}
The problem of counting the number of populations for a special
choice of relative shifts $\bs b^{(i)}$ can be approached
via Schubert calculus.

Let $\mc V$ be a complex vector space of dimension $d+1$  and
\be
\mathcal F=\{ 0\subset F_1\subset F_2\subset\dots\subset
F_{d+1} =\mc V\}, \qquad \dim F_i=i, 
\ee
 a full flag in $\mc V$.  Let $Gr(N+1,\mc V)$ be the Grassmanian
of all $(N+1)$-dimensional subspaces in $\mc V$.

Let $\bs a =(a_1,\dots,a_{N+1})$, 
$ d - N \geq a_1\geq a_2\geq \dots \geq a_{N+1}\geq 0$,
be a non-increasing sequence of  non-negative integers.
Define the  {\it Schubert cell} $G^0_{\bs a}(\mc F)$
associated to the flag $\mathcal F$ and sequence $\bs a$ as the set
\bea
 \{V\in Gr(N+1,\mc V)\ & | &
\dim (V \cap F_{d- N + i-a_i}) = i ,
\\
&& \dim (V \cap F_{d- N + i-a_i-1}) = i-1, 
\ \text{for}\ 
i = 1, \dots , N+1 \} .
\eea
The closure $G_{\bs a}(\mc F)$ of the Schubert cell is called {\it
the Schubert cycle}. For a fixed flag $F$, the Schubert cells 
form a cell decomposition of the Grassmanian.
The codimension of $G^0_{\bs a}(\mc F(z))$
in the Grassmanian is $|\bs a| = a_1 + \dots + a_{N+1}$. 

Let $\mc V = \C_d[x]$ be the space of polynomials
of degree not greater than $d$, dim $\mc V = d + 1$. For any  
$z \in \C \cup \infty$,
define a full flag in $\C_d[x]$,
\be
\mathcal F(z)\ = \ \{ 0 \subset F_1(z) \subset F_2(z) \subset \dots \subset
F_{d+1}(z)=\mc V \}\ .
\ee 
For $z \in \C $ and any $i$, we set $F_i (z)$ to be the subspace of all
polynomials divisible by $\prod_{j=1}^{d+1-i}( x - z-(j-1)h)$. 
For any $i$, we set $F_i(\infty) $ to be the subspace of all
polynomials of degree less than $i$.

Let $V \in Gr (N+1, \C_d[x])$.
For any $z \in \C \cup \infty$, let $\bs a(z)$ be such a unique
sequence  that $V$ belongs to the cell $G^0_{\bs a (z)}(\mc F(z))$.
We say that a point $z \in \C \cup \infty$ is an {\it $h$-ramification
  point } for $V$, if  $\bs a(z) \neq (0, \dots , 0)$. 
We call  $\bs a(z)$ the {\it ramification condition of $V$ at $z$}.

If $u_1,\dots, u_{N+1}$ is a basis of $V$ 
then the ramification points are zeroes of the discrete Wronskian
$W(x)=W(u_1,\dots,u_{N+1})(x)$. Indeed, given $z\in\C$,
without loss of generality we can assume that
$u_i\in\F_{d+2-i-a_i(z)}(z)$ and therefore the matrix $\{u_i(x+(j-1))\}_{i,j=1}^{N+1}$ is upper triangular. Moreover, this matrix has 
a zero diagonal element if and only if $|a|>0$.

Fix $h$-ramification conditions at 
$z_1, \dots , z_n, \infty$ so that
\be\label{ram cond}
\sum _{s=1}^n | \bs a (z_s)| \ + \ | \bs a (\infty)|\ =\ 
\dim \ Gr ( N+1, \C_d[x]) \ .
\ee

\medskip

\noindent {\bf Counting Problem.} 
{\it Compute the number of spaces of polynomials of dimension $N+1$ 
with these ramification properties.}

\medskip

In other words we ask to count the number of points in the
intersection of Schubert cycles. The intersection index of the
Schubert cycles is well known.

Namely, for a non-increasing sequence $\bs a =(a_1,\dots,a_{N+1})$, 
$  a_1\geq a_2\geq \dots,\geq a_{N+1}\geq 0$,
of  non-negative integers, denote
$\tilde L_{\bs a}$ the finite dimensional irreducible
$gl_{N+1}$-module with highest weight $\bs a$. Let $L_{\bs a}$ be the
$sl_{N+1}$ module obtained by restriction of $\tilde L_{\bs a}$.

\begin{theorem}\label{multiplicity}{\rm (\cite{F})}
The intersection index of Schubert cycles $G_{\bs a(z_s)}(F(z_s))$,
$s=1,\dots,n$, and 
$G_{\bs a(\infty)}(F(\infty))$ equals
the multiplicity 
of the trivial $sl_{N+1}$-module in the tensor product of $sl_{N+1}$-modules
\be
 L_{\bs a (z_1)} \otimes \dots  \otimes
L_{\bs a (z_n)} \otimes  L_{\bs a (\infty)} \ .
\ee
\end{theorem}

Conjecturally, for almost all $z_1, \dots z_n$ the number of spaces
$V$ with such ramification conditions is equal to the above multiplicity. 

Fix an initial data $\bs z, \bs \La, \bs b^{(i)}$.
Until the end of 
this section we assume that $z_s-z_r\not \in h\Z$ for all $s\neq r$ and
\bean\label{sl shift}
b_s^{(j)}=-\sum_{i=1}^j\La_s^{(i)}.
\eean
Fix an $sl_N$ weight and write it in the form $w\cdot \La_\infty$ where
$\La_\infty$ is dominant integral and $w$ is an element of Weyl
group. 

Let $\bs y$ be an $h$-critical point associated to the initial data
and the weight at infinity $w\cdot \La_\infty$. 
Let $P$ be the population of $h$-critical
points originated at $\bs y$ and let $V_P$ be the corresponding
fundamental space. Let $d$ be large enough, so that $V_P\subset\C_d[x]$.

\begin{theorem}
The points $z_1,\dots,z_n$ and
$\infty$ are ramification points of $V_P$. 
The ramification condition at $z_s$ is
$a_i(z_s)=\sum_{j=1}^{N+1-i}\La_s^{(j)}$. The ramification condition at
$\infty$ is $a_i(\infty)=d-N-l_1-\sum_{j=1}^{i-1}\La_\infty^{(j)}$,
where $l_1$ is defined from $\bs \La$ and $\La_\infty$ via \Ref{w at
  inf}.
\end{theorem}
\begin{proof}
Follows from Lemma \ref{wr(u) lem}, cf. proof of Lemmas 5.8, 5.10 
in \cite{MV}.
\end{proof}

For an integral dominant $sl_{N+1}$ weight $\La$, we denote $L_\La$
the irreducible $sl_{N+1}$ module with highest weight $\La$.

\begin{cor}. The number of (discretely lying) 
populations associated to initial data $\bs z_i,\bs \La, \bs b^{(i)}$
such that \Ref{sl shift} holds, is not greater
than the multiplicity of $L_{\La_\infty}$ in the tensor product
$L_{\La_1}\otimes\dots\otimes L_{\La_n}$.
\end{cor}
\begin{proof}
The corollary holds since the number of isolated points of the
intersection of the corresponding Schubert cycles is not greater than the
intersection index of the cycles.
\end{proof} 

\subsection{Generating morphism}\label{flag section}
We identify a population with the variety of full flags in the
fundamental space.

Let $V = V_P$ be the fundamental space of a  population $P$. By
Proposition \ref{T form frame}, the 
polynomials $T_i$ defined in \Ref{T polyn} form a frame of $V_P$.

Let $FL(V)$ be the variety of all full flags
\be
\mathcal F \ = \ \{ 0 \subset F_1 \subset F_2 \subset \dots \subset
F_{N+1} = V \}\ 
\ee 
in $V$. For any $\mathcal F \in FL(V)$ define an $N$-tuple of polynomials
$\bs y^{\mathcal F} =  (y_1^{\mathcal F}, \dots , y_N^{\mathcal F})$\ 
as follows. Let $u_1, \dots , u_{N+1}$ be a basis in $V$ such that for
any $i$ the polynomials 
$u_1, \dots , u_i$ form a basis in $F_i$. 
We say that this basis {\it is adjusted } to the flag $\mc F$ and the flag
$\mc F$ {\it is generated} by the basis $u_1, \dots , u_{N+1}$. 

Define the polynomials 
\bean\label{div}
y_i^{\mathcal F}(x) \ = \ \frac{ W(u_1, \dots , u_i) (x)}
{ \prod _{j=1}^{i-1}T_1(x+(j-1)h)\ 
\prod_{j=1}^{i-2} T_2(x+(j-1)h)\  \dots\  T_{i-1}(x)} \ . 
\eean

The correspondence $\mc F \mapsto \bs y^{\mc F}$ gives a morphism
\bean\label{morphism}
\beta \ : \ FL(V) \ \to \PCN , 
\eean
called {\it the generating morphism of $V$}.

\begin{theorem}\label{Main}
Let $P$ be a population of $h$-critical points with the fundamental space $V$.
Then the generating morphism  defines an isomorphism of $FL(V)$ and
the population $P\subset\PCN$.
\end{theorem}
\begin{proof}
Proof is the same as the proof of the first part of Theorem 5.12 in \cite{MV}.
\end{proof}

\medskip

{\bf Remark.}
Recall that $V$ is the $(N+1)$-dimensional complex vector space of
polynomials which is contained in the kernel of the fundamental
operator $D_P$ and $D_P$ is a linear difference operator of order $N+1$.
Therefore the full flags in $V$ also label the decompositions of 
$D_P$ to $N+1$ linear factors of the form $(\partial -
f(x))$, where $f(x)$ is a rational function of $x$. Thus, $h$-critical
points are in bijective correspondence with such factorizations of
the fundamental operator $D_P$. 

\medskip

Fix a flag $\mc F^0\in FL(V)$.
For any $\mc F \in FL(V)$ define a permutation
$w(\mc F)$ in the symmetric group $S^{N+1}$ as follows. 
Define $w_1(\mc F)$ as the minimum of $i$ 
such that  $F_1 \subset F^0_i$. Fix a basis vector $u_1 \in F_1$. 
Define $w_2(\mc F)$ as the minimum of $i$ such that  
there is a basis in $F_2$ of the form $u_1, u_2$ with $u_2 \in F_i^0$.
Assume that 
$w_1(\mc F), \dots , w_j(\mc F)$ and $u_1,\dots,u_j$ are determined.
Define $w_{j+1}(\mc F)$ as the minimum of $i$ such that  
there is a basis in $F_{j+1}$ of the form $u_1, \dots , u_j, u_{j+1}$ with 
$u_{j+1} \in F_i^0$. As a result of this 
procedure we define $w(\mc F) = (w_1(\mc F), \dots , w_{N+1}(\mc F))
\in S^{N+1}$ 
and a basis $u_1, \dots , u_{N+1}$ which generates $\mc F$ and such that
$u_i \in F^0_{w_i(\mc F)}$.

For $w \in S^{N+1}$, define 
\be
G^{\mc F^0}_w  = \{ \mc F \in FL(V), \ w(\mc F) = w \} .
\ee 
The algebraic variety $G^{\mc F^0}_w$ is
called {\it the Bruhat cell } associated with $\mc F^0$ and  $w \in S^{N+1}$.
The set of all Bruhat cells form a cell decomposition of $FL(V)$:
\be
FL(V)=\sqcup_{w\in S^{N+1}}G^{\mc F^0}_w.
\ee

Recall that the symmetric group $S^{N+1}$ is identified with the
$sl_{N+1}$ Weyl group in such a  way that the simple transposition
$(i,i+1)$ corresponds to the simple reflection with respect to
$\al_i$.

Now we are ready to describe the set of $(N+1)$-tuples in a population of
a  fixed degree.
Let $\bs y$ represent an $h$-critical point associated with the
initial data and the weight at infinity $w\cdot \La_\infty$, where 
$\La_\infty$ is dominant integral. Let $P=P(\bs y)$ be 
the corresponding population and $\beta:FL(V)\to P$ the generating
isomorphism. 

\begin{theorem}
The set of $N$-tuples $\bar{\bs y}$ in $P(\bs y)$ associated with a
weight at infinity $w\cdot\La_\infty$, where $\La_\infty$ is integral
dominant, coincides with the image of the Bruhat cell $\beta(G_w^{\mc
  F(\infty)})$. 
\end{theorem}
\begin{proof}
The theorem is proved similar to Corollary 5.23 in \cite{MV}.
\end{proof}

In particular, each population contains at most one tuple $\bs y$
which represents an $h$-critical point associated to a integral
dominant weight at infinity. Thus the number of critical points
associated to an initial data and an integral dominant weight is
bounded from above by the number of the populations. We conjecture
that for
generic values of $z_i$ this bound is exact.

\section{$h$-selfdual vector spaces of polynomials}\label{selfdual sec}

\subsection{Dual spaces}
Let $V$ be a space of polynomials of dimension $N+1$ 
with frame $T_1(x), \dots , T_N(x)$.
For $u_1, \dots , u_i \in V$, the polynomial
\be
W^{\dagger}(u_1, \dots , u_i)(x)\ := \ \frac{W(u_1, \dots , u_i)(x)}
{\prod _{j=1}^{i-1}T_1(x+(j-1)h)\ \prod_{j=1}^{i-2} T_2(x+(j-1)h)\
  \dots\  T_{i-1}(x)} 
\ee
is called {\it the divided Wronskian}. 

Note that if polynomials 
$u_1,\dots, u_{N+1}$ form a basis of $V$ then the divided 
Wronskian $W^\dagger(u_1,\dots,u_{N+1})$ is a non-zero constant.

Let $V^\dagger$ be the set of polynomials
$W^\dagger(u_1,\dots,u_N)$, where $u_i\in V$. Clearly $V^\dagger$ is a
vector space of dimension $N+1$. 
We call $V^\dagger $ the {\it $h$-dual space of V}.

We have a non-degenerate pairing
\be
V \otimes V^\dagger  \to \C,
\qquad u \otimes W^\dagger (u_1, \dots , u_N) 
\mapsto  W^\dagger(u, u_1, \dots , u_{N})(x).
\ee

For $i = 1, \dots , N,$ set 
\bean\label{tilde T}
T^\dagger_{i}(x) = T_{N+1-i}(x + (i-1)h).
\eean

The following two lemmas are obtained 
from Wronskian identities in Lemmas
\ref{wr id 1} and \ref{f wronskian}.

\begin{lem}
The polynomials 
$T^\dagger_1(x), \dots , T^\dagger_N(x)$ form a frame of $V^\dagger$.
\hfill $\square$
\end{lem}

For a space of polynomials $V$ and $a\in\C$, we denote $V(x+a)$ the space of
polynomials spanned by
$u(x+a)$, $u(x)\in V$.

\begin{lem}
We have $V^{\dagger \dagger} = V(x+(N-1)h)$.
\hfill $\square$
\end{lem}

\subsection{$h$-selfdual spaces and canonical bilinear form}
We say that $V$ is {\it $h$-selfdual}, if $V^\dagger = V(x-(N-1)h/2)$. 

A space $V$ is $h$-selfdual if and only if $V^\dagger$ is $h$-selfdual.

If $V$ is $h$-selfdual, then 
\bean\label{T=tilde T}
T_i(x) \ =\ T_{N+1-i}(x + (i-1)h - (N-1)h/2). 
\eean
In particular, if $T_i$ are of the form \Ref{T polyn}, then 
\be
\La_s^{(i)}=\La_s^{(N+1-i)},\qquad b_s^{(i)}=b_s^{(N+1-i)}+(i-1)h-(N-1)h/2.
\ee

For instance, the space of polynomials of degree not greater than $N$ is
$h$-selfdual. In this case all polynomials $T_i$ are equal to $1$.

Let $V$ be $h$-selfdual. 
Define a non-degenerate pairing $( \ ,\  ) :  V \otimes V  \to \C$. 
If $u, v \in V$, then we write 
$v(x) = W^{\dagger} (u_1, \dots , u_{N})(x - (N-1)h/2)$ with $u_i\in V$
 and set 
$( u , v ) = W^{\dagger}(u, u_1 , \dots , u_{N})$.
This pairing is called the {\it canonical bilinear form}.

A basis $u_1, \dots , u_{N+1}$ in a space of polynomials $V$ is called {\it
a  Witt basis} if for $i = 1, \dots , N+1$ we have
\bean\label{dar 2 eqn}
u_i (x)\ =\ W^\dagger (u_1, \dots , \widehat{u}_{N+2-i}, \dots ,
u_{N+1})(x - (N-1) h/2).
\eean
Clearly, if $V$ has a Witt basis, then $V$ is $h$-selfdual.

\begin{theorem}\label{special basis prop}
Let $V$ be $h$-selfdual. Then $V$ has a Witt basis. In particular, 
the form $(\ ,\ ) :  V \otimes V  \to  \C$
is symmetric if the dimension of $V$ is odd 
and skew-symmetric if the dimension of $V$ is even.
\end{theorem}
\begin{proof}
The theorem is proved similarly to Theorem 6.4 in \cite{MV}.
\end{proof}

\subsection{Isotropic flags}\label{isotropic}
Let $u_1, \dots , u_{N+1}$ be a basis in a vector space
$V$ of polynomials. Denote 
$W_i(x) = W^\dagger (u_1, \dots , \widehat{u}_{i}, 
\dots , u_{N+1})(x - (N-1)h/2)$. 
\begin{lem}\label{suf self}
If for $ i = 1, \dots , N$
\bean\label{dar 1 eqn}
W^\dagger (u_1, \dots , u_i)(x)\ =\ 
c_i\ W^\dagger (u_1, \dots , u_{N+1-i})(x + (i-1) h - \frac{N-1}{2} h), 
\eean
where $c_1, \dots , c_N$ are some non-zero complex numbers, 
then for every $i$, the polynomial
$u_i(x)$ is a linear combination of $W_{N+1}(x),  W_N(x),
 \dots , W_{N+2-i}(x)$.
\end{lem}
\begin{proof} 
The proof of the lemma is similar to the proof of Theorem 6.8 in \cite{MV}.
\end{proof}
\begin{cor} \label{slf cor}
If $V$ has a basis satisfying \Ref{dar 1 eqn}, then $V$ is $h$-selfdual.
\end{cor}

Let $V$ be an $h$-selfdual space of polynomials. 
For a subspace $U \subset V$ denote $U^\perp$ its  orthogonal
complement.
A full flag $\mathcal F = \{0\subset F_1 \subset \dots \subset F_{N+1} = V\}$ 
is called {\it isotropic}  if $F_i^{\perp} = F_{N+1-i}$ for  
$i = 1, \dots , N$. 

\begin{prop}
Let $\mc F$ be a flag in an $h$-selfdual space of polynomials $V$ 
and $u_1, \dots , u_{N+1}$  a basis in $V$ adjusted to $\mc F$. 
Then $\mc F$ is isotropic if and only if \Ref{dar 1 eqn} holds. 
\end{prop}
\begin{proof}
If $F_i^\perp = F_{N+1-i}$ then we have two bases in $F_{N+1-i}$: 
the basis $u_1, \dots , u_{N+1-i}$ and the basis
$W_{N+1}, W_N, \dots , W_{i+1}$. Hence
\bea
W^\dagger (u_1, \dots , u_{N+1-i})(x)  &=& 
 \text{const} W^\dagger (W_{N+1}, W_N, \dots , W_{i+1})(x)=
\\
&=& \text{const} W^\dagger (u_1,  \dots , u_i)(x + (N-i) h
 - (N-1)h/2).
\eea
The only if part of the lemma is proved.

Now, let \Ref{dar 1 eqn} hold. We prove that the spaces spanned by
$u_1(x),\dots,u_i(x)$ and by 
$W_{N+1}(x-(N-1)h/2),\dots,W_{N-i}(x-(N-1+2)h/2)$ 
are the same by induction on $i$. 
For $i=1$ it is just equation \Ref{dar 1 eqn}. Assume that the
statement is proved for $i=1,\dots,i_0-1$. Then
\bea
W^\dagger (u_1, \dots , u_{i_0})(x)={\rm const}  W^\dagger (u_1,
\dots , u_{N+1-i_0})(x + (i_0-1) h - (N-1)h/2)=\\
{\rm const}
 W^\dagger (W_{N+1}, W_N, \dots , W_{N-i_0+2})(x-(N-1)h/2)=\\
{\rm const} 
 W^\dagger (u_1(x),\dots,u_{i_0-1}(x),W_{N-i_0+2}(x-(N-1)h/2)),
 \eea
and the step of induction follows.
\end{proof}

\section{The case of $B_N$}\label{B sec}
Consider the root system of type $B_N$ corresponding to the Lie algebra
$so_{2N+1}$.
Let $\al_1,\dots,\al_{N-1}$ be long simple roots and $\al_N$ the short
one. We have
\be
(\al_N,\al_N)=2,\qquad (\al_i,\al_i)=4,\qquad (\al_i,\al_{i+1})=-2,
\ee
for $i=1,\dots,N-1$, and all other scalar products are zero. 

\subsection{Definition of critical points of $B_N$ type}\label{def B section} 
We fix 
ramification points $\bs z=(z_1,\dots,z_n)\in\C^n$, non-zero
dominant integral $B_N$
weights $\bs \La=(\La_1,\dots,\La_n)$ and relative shifts $\bs
b^{(i)}=(b^{(i)}_1,\dots,b^{(i)}_n)\in\C^n$, $i=1,\dots,N$. 
We call this data a $B_N$ initial data. We also set
$\La_s^{(i)}=(\La_s,\al_i^\vee)$. 

Given $\bs l=(l_1,\dots,l_N)\in\Z_{\geq 0}^N$, we define the $B_N$ 
weight at infinity $\La_\infty$ by the formula \Ref{w at inf}.

For a $B_N$ initial data $\bs z,\bs \La, \bs b^{(i)}$ 
and a weight at infinity 
$\La_\infty$, we define an $sl_{2N}$ initial data which
consists of ramification points $z_i$, non-zero dominant integral $sl_{2N}$
weights $\La_i^A$,
relative shifts $\bs b^{(i),A}$ and we also define an $sl_{2N}$ 
weight at infinity $\La_\infty^A$. 
Namely, given a $B_N$ weight $\La$, the $sl_{2N}$ weight $\La^A$ is defined by
\be
(\La^A,\al_i^A)=(\La^A,\al_{2N-i}^A)=(\La,\al_i^\vee),
\ee
where $i=1,\dots,N$ and $\al_i^A$ are roots of $sl_{2N}$. The shifts
$\bs b^{(i),A}$ are defined by
\be
b_s^{(i),A}=b_s^{(2N-i),A}+(i-N)h=b_s^{(i)}.
\ee

Given a set of complex numbers $t_j^{(i)}$, $i=1,\dots,N$,
$j=1,\dots l_i$, we represent it by the $N$-tuple of polynomials $\bs
y=(y_1,\dots, y_N)$, where $y_i(x)=\prod_{j=1}^{l_i}
(x-t_j^{(i)})$. We define the corresponding $(2N-1)$-tuple $\bs y^A$
by the formula
\bean\label{A=B}
y_i^A(x)=y_i(x), \qquad y_N^A(x)=y_N(x),\qquad
y_{i+N}^A(x)=y_{N-i}(x+ih),
\eean
where $i=1,\dots,N-1$.

We propose the following definition of the $h$-critical points of type
$B_N$, cf. Lemma 7.1 in \cite{MV}. 
The set of complex numbers $t_j^{(i)}$, $i=1,\dots,N$,
$j=1,\dots, l_i$, is called an {\it $h$-critical point of $B_N$ type}
associated to the initial data $\bs z$, $\bs \La$, $\bs b^{(i)}$
and the weight at infinity $\La_\infty$ if the
corresponding $(2N-1)$-tuple $\bs y^A$ represents an $sl_{2N}$
$h$-critical point associated to the initial data $\bs z$, $\bs
\La^A$, $\bs b^{(i),A}$ and the weight at infinity $\La_\infty^A$.

In this case we 
say that the $N$-tuple $\bs y$ represents an $h$-critical point of
the $B_N$ type.

Equivalently, the set of numbers 
$t_j^{(i)},$ $i=1,\dots,N$,
$j=1,\dots, l_i$, is an
$h$-critical point of $B_N$ type if it satisfies the following system
of algebraic equations:
\begin{align*}
\prod_{s=1}^n \frac {t_j^{(i)}-z_s + b_s^{(i)} h + \La_s^{(i)} h}
{t_j^{(i)}-z_s + b_s^{(i)} h} \
\prod_{k=1}^{l_{i-1}}
\frac{t_j^{(i)} - t_k^{(i-1)} + h}
{t_j^{(i)} - t_k^{(i-1)}} \
& \prod_{k \neq j}
\frac {t_j^{(i)} - t_k^{(i)} - h}
{t_j^{(i)} - t_k^{(i)} + h} \times  \\
& \times \prod_{k=1}^{l_{i+1}}
\frac{t_j^{(i)} - t_k^{(i+1)}}
{t_j^{(i)} - t_k^{(i+1)} - h} =  1, \\
\prod_{s=1}^n \frac {t_j^{(N)}-z_s + b_s^{(N)} h + \La_s^{(N)} h}
{t_j^{(N)}-z_s + b_s^{(N)} h} \
\prod_{k=1}^{l_{N-1}}
\left(\frac{t_j^{(N)} - t_k^{(N-1)} + h}
{t_j^{(N)} - t_k^{(N-1)}}\right)^2 
& \prod_{k \neq j}
\frac {t_j^{(N)} - t_k^{(N)} - h}
{t_j^{(N)} - t_k^{(N)} + h}=1,
\end{align*}
where $i=1,\dots,N-1$.

Note that in the quasiclassical limit $h\to 0$, this system 
becomes system (2.2) of \cite{MV} specialized to the case of 
$B_N$. 

\subsection{Reproduction procedure of $B_N$ type}
We describe the concepts of reproduction and population. 
All of that follows from the definitions in the case of $sl_{2N}$.

An $N$-tuple $\bs y$ is called fertile in $B_N$ sense is there exist 
polynomials
$\tilde y_i$, $i=1,\dots,N-1$ and $\tilde y_N$ such that
\bea
W(y_i,\tilde y_i)(x)&=&y_{i+1}(x)y_{i-1}(x+h)T_i(x),\\
W(y_N,\tilde y_N)(x)&=&y_{N-1}^2(x+h)T_N(x).
\eea
The $N$-tuple $\bs y^{(i)}=(y_1,\dots,\tilde y_i,\dots,y_N)$ is called
an immediate descendant of $\bs y$ in the direction $i$. Immediate
descendants $\bs y^{(i)}$ are also fertile in $B_N$ sense.

From Lemma \ref{fertile=critical} we obtain that an $N$-tuple $\bs y$
represents a critical point if and only if it is fertile and $\bs y^A$
is generic. In this case if $(\bs y^{(i)})^A$ is generic, then $\bs y^{(i)}$ 
also represents an $h$-critical
point of type $B_N$.

The $B_N$ population originated at a critical point $\bs y$ is the
minimal set $P(\bs y)$ of fertile $N$-tuples such that $\bs y \in P(\bs y)$ 
and if $\bar{\bs y}\in P(\bs y)$ then the immediate descendants 
${\bar{\bs y}}^{(i)}\in P(\bs y)$. 

Obviously, the $B_N$ population is contained in the corresponding
$sl_{2N}$ population: if $\bar{\bs y}\in P(\bs y)$ then  
$\bar{\bs y}^A\in P_A(\bs y^A)$, where we denote $P_A(\bs y^A)$ the
$sl_{2N}$ population of critical points originated at $\bs y^A$.
The fundamental space of
the $sl_{2N}$ population $P_A(\bs y^A)$ contains a flag with the
property \Ref{dar 2 eqn} and therefore is $h$-selfdual. We call this
space the {\it fundamental space of the $B_N$ population $P(\bs y)$}.

The corresponding fundamental operator is the difference operator of
order $2N$ obtained from \Ref{Diff} by substituting \Ref{A=B} and
\Ref{T=tilde T}:
\bea
\prod^{N \to 1}_i
\left( \partial  -  
\frac{y_{N-i}(x+(i+1)h)}{y_{N-i}(x+ih)} 
\frac{y_{N+1-i}(x+(i-1)h)}{y_{N+1-i}(x+ih)}
T(x+ih)\prod_{s=1}^{i-1} \frac{T_{N-s}(x + ih)}{T_{N-s}(x+(i-1)h)} 
 \right)\times\\
\prod^{1 \to N}_i
\left( \partial  -  
\frac{y_{N+1-i}(x+h)}{y_{N+1-i}(x)} 
\frac{y_{N-i}(x)}{y_{N-i}(x+h)}
\prod_{s=1}^{N-i} \frac{T_{s}(x + (N-i-s+1)h)}{T_{s}(x+(N-i-s)h)}
 \right),
\eea
where $T(x)=\prod_{s=1}^{N} T_{s}(x + (N-s)h)/T_{s}(x+(N-s-1)h)$.

Clearly, two populations of $h$-critical points of type $B_N$ either
do not intersect or coincide. 

\subsection{A $B_N$ population and the $C_N$ flag variety}
Let $\bs y$ be an $h$-critical point of $B_N$ type and let $V$ be the
fundamental space of the population $P_A(\bs y^A)$. Then $V$ is an
$h$-selfdual space of dimension $2N$ and therefore has a non-degenerate 
skew-symmetric canonical bilinear form.

The special symplectic Lie algebra of $V$ consists of all
traceless endomorphisms $x$ of $V$
such that $(xv , v') + (v , xv') = 0$ for all $v, v' \in V$.
Let $\bs u = (u_1, \dots , u_{2N})$ be a Witt basis in $V$. We have
$( u_i, u_{2N+1-i}) = (-1)^{i+1}$, $i = 1, \dots , N$ ,
and $(u_i , u_j) = 0$ if $i + j \neq 2N+1$.
This choice of basis identifies the special symplectic Lie algebra 
with a Lie subalgebra of $sl_{2N}$, which is denoted
$sp_{2N}$. The Lie algebra $sp_{2N}$ has the root system of type $C_N$.

Denote $E_{i, j}$ the matrix with zero entries except 1 at the
intersection of the $i$-th row and the $j$-th column.
The lower triangular part of $sp_{2N}$ is spanned by
matrices $ E_{i, i+1} + E_{ 2N-i, 2N+1-i}$  for $i = 1, \dots, N - 1$
and $E_{N, N+1}$. Denote these matrices $X_1, \dots , X_N$, respectively.

Now we describe the action of one-parametric subgroups in the special
symplectic group in the 
direction of $X_i$ on the basis $\bs u=(u_1,\dots,u_{2N})$. 
We have
\bea
&& e^{c X_i} \bs u =\\
&&=(u_1, \dots , u_{i-1}, u_i + c u_{i+1}, u_{i+1}, \dots ,
u_{2N-i}, u_{2N+1-i} + c u_{2N+2-i}, u_{2N+2-i}, \dots , u_{2N}),\\
&&e^{c X_i} \bs u  = 
(u_1,  \dots , u_{N-1}, u_N + c u_{N+1}, u_{N+1},
  \dots , u_{2N}),
\eea
where $i=1,\dots,N-1$. 

We also set
\bea
e^{\infty X_i} \bs u &=&
(u_1, \dots , u_{i-1}, u_{i+1}, u_{i}, \dots ,
u_{2N-i}, u_{2N+2-i}, u_{2N+1-i}, \dots , u_{2N}),\\
e^{\infty X_i} \bs u & = &
(u_1,  \dots , u_{N-1}, u_{N+1}, u_{N},
  \dots , u_{2N}).
\eea

Note $e^{cX_i}\bs u$ is a Witt basis
for all $i,c$ and the corresponding flag is isotropic.

Given a basis $\bs u=(u_1,\dots,u_{2N})$ of $V$ we denote 
$\mc F(\bs u)$ the full flag generated by this basis.
Let $FL(V)$ be the variety of all full flags of $V$.

Denote $FL^\perp(V) \subset FL(V)$ the subvariety of all isotropic
flags and define the  $B_N$ generating morphism  
$\beta : FL^\perp(V) \to  \PCN$, $\mathcal F(\bs u) \mapsto
(y_1^{\mc F},\dots,y_N^{\mc F})$, where 
$y_i^{\mc F}=W^\dagger(u_1,\dots,u_i)$.

\begin{lem}\label{desc lem}
If $F(\bs u)$ is an isotropic flag then $\beta(\mc F(\bs u))$ is
fertile. Moreover the set of all
immediate descendants of $\beta(\mc F(\bs u))$ in the direction $i$ coincides
with the set $\beta (\mc F(e^{cX_i}\bs u))$, $c\in\bar{\C}$.
\end{lem}
\begin{proof}
Follows from the Wronskian identities
\bea
W(y_i^{\mc
  F},W^\dagger(u_1,\dots,u_{i-1},u_{i+1}))(x)&=&T_i(x)y_{i-1}^{\mc
  F}(x+h)y_{i+1}^{\mc F}(x),\\
W(y_N^{\mc
  F},W^\dagger(u_1,\dots,u_{N-1},u_{N+1}))(x)&=&T_N(x)(y_{N-1}^{\mc
  F}(x+h))^2,
\eea
where $i=1,\dots, N-1$.
\end{proof}

\begin{theorem}
The generating morphism $\beta$ of type $B_N$ gives rise to the isomorphism
$FL^\perp(V)\to P(\bs y)\subset\PCN$.
\end{theorem}
\begin{proof}
The generating morphism $\beta$ is an isomorphism of $FL^\perp(V)$ to the
image by Theorem  \ref{Main}. Clearly, the image of $\beta$ contains the
population $P(\bs y)$. The image of $\beta$ coincides with $P(\bs y)$ by
Lemma \ref{desc lem}.
\end{proof} 

Note that $FL^\perp$ is isomorphic to the flag variety of the special
symplectic group which corresponds to the $C_N$ root system. This flag
variety has a Bruhat cell decomposition
$FL^\perp(V)=\sqcup_{w\in \mc W}G_w^{\mc F(\infty),C}$, 
where $\mc W$ is the $C_N$ Weyl group.
Note that the Weyl groups of $B_N$ and $C_N$ are canonically identified.
See for details, e.g. \cite{MV}, section 7.3.

As before, the set of all $N$-tuples in a $B_N$ population of the same degree
is the Bruhat cell in the $C_N$ flag variety.
\begin{theorem}
The set of $N$-tuples $\bar{\bs y}$ in $P(\bs y)$ associated with a
weight at infinity $w\cdot\La_\infty$, where $\La_\infty$ is integral
dominant, coincides with the image of the Bruhat cell $\beta(G_w^{\mc
  F(\infty),C})$. 
\end{theorem}
\begin{proof}
Completely parallel to the proof of Corollary 7.12 in \cite{MV}.
\end{proof}

\section{The case of $C_N$}\label{C sec}
Consider the root system of type $C_N$ corresponding to the Lie algebra
$sp_{2N}$.
Let $\al_1,\dots,\al_{N-1}$ be short simple roots and $\al_N$ the long
one. We have
\be
(\al_N,\al_N)=4,\qquad (\al_i,\al_i)=2,\qquad (\al_{N-1},\al_N)=-2 
\qquad (\al_{i-1},\al_{i})=-1,
\ee
for $i=1,\dots,N-1$ and all other scalar products are zero. 

\subsection{Definition of critical points of $C_N$ type}\label{def C section} 
We fix ramification points $\bs z=(z_1,\dots,z_n)\in\C$, 
dominant integral $C_N$
weights $\bs \La=(\La_1,\dots,\La_n)$ and relative shifts $\bs
b^{(i)}=(b_1^{(i)}, \dots,\bs b_n^{(i)})\in\C^n$, 
$i=1,\dots,N$. We call this data a $C_N$ initial data. We also set
$\La_s^{(i)}=(\La_s,\al_i^\vee)$. 

Given $\bs l=(l_1,\dots,l_N)\in\Z_{\geq 0}^N$, we define the $C_N$ 
weight at infinity $\La_\infty$ by the formula \Ref{w at inf}.

For a $C_N$ initial data $\bs z,\bs \La,\bs b^{(i)}$ and a weight at
infinity $\La_\infty$, 
we define an $sl_{2N+1}$ initial data which
consists of ramification points $z_i$, dominant integral $\La_i^A$,
relative shifts $\bs b^{(i),A}$ and we also define an $sl_{2N+1}$
weight at infinity $\La_\infty^A$. Namely,
given a $C_N$ weight $\La$, the $sl_{2N+1}$ weight $\La^A$ is defined by
\be
(\La^A,\al_i^A)=(\La^A,\al_{2N+1-i}^A)=(\La,\al_i^\vee),
\ee
where $i=1,\dots,N$ and $\al_i^A$ are roots of $sl_{2N+1}$. The shifts
$\bs b^{(i),A}$ are defined by
\be
b_s^{(i),A}=b_s^{(2N+1-i),A}+(i-N+1/2)h=b_s^{(i)}.
\ee

Given a set of complex numbers $t_j^{(i)}$, $t_k^{(N)}$ where 
$i=1,\dots,N-1$,
$j=1,\dots, l_i$, $k=1,\dots,2l_N$, 
we represent it by the $N$-tuple of polynomials $\bs
y=(y_1,\dots, y_N)$, where $y_i(x)=\prod_{j=1}^{l_i}
(x-t_j^{(i)})$, $y_N=\prod_{j=1}^{2l_N}
(x-t_j^{(N)})$. We define the corresponding $(2N)$-tuple $\bs y^A$
by the formula
\bean\label{A=C}
y_i^A(x)=y_i(x), \qquad 
y_{i+N}^A(x)=y_{N+1-i}(x+(i-1/2)h),
\eean
where $i=1,\dots,N$.

We propose the following definition of the $h$-critical points of type
$C_N$, cf. Section 7.2 in \cite{MV}. 
The set of complex numbers $t_j^{(i)}$, $t_k^{(N)}$, where $i=1,\dots,N-1$,
$j=1,\dots l_i$, $k=1,\dots,2l_N$,
is called an {\it $h$-critical point of $C_N$ type}
associated to the initial data $\bs z$, $\bs \La$, $\bs b^{(i)}$
and the weight at infinity $\La_\infty$ if the
corresponding $(2N)$-tuple $\bs y^A$ represents an $sl_{2N+1}$
$h$-critical point associated to the initial data $\bs z$, $\bs
\La^A$, $\bs b^{(i),A}$ and the weight at infinity $\La^A_\infty$.

In this case we 
say that the $N$-tuple $\bs y$ represents an $h$-critical point of the
$C_N$ type.

Equivalently, the set of numbers 
$t_j^{(i)}$, $t_k^{(N)}$, where $i=1,\dots,N-1$,
$j=1,\dots, l_i$, $k=1,\dots,2l_N$, is an
$h$-critical point of $C_N$ type if it satisfies the following system
of algebraic equations:
\begin{align*}
\prod_{s=1}^n \frac {t_j^{(i)}-z_s + b_s^{(i)} h + \La_s^{(i)} h}
{t_j^{(i)}-z_s + b_s^{(i)} h} \
\prod_k
\frac{t_j^{(i)} - t_k^{(i-1)} + h}
{t_j^{(i)} - t_k^{(i-1)}} & 
 \prod_{k \neq j}
\frac {t_j^{(i)} - t_k^{(i)} - h}
{t_j^{(i)} - t_k^{(i)} + h} \times  \\
& \times \prod_k
\frac{t_j^{(i)} - t_k^{(i+1)}}
{t_j^{(i)} - t_k^{(i+1)} - h} =  1, \\
\prod_{s=1}^n \frac {t_j^{(N)}-z_s + b_s^{(N)} h + \La_s^{(N)} h}
{t_j^{(N)}-z_s + b_s^{(N)} h} \
\prod_{k=1}^{l_{N-1}}
\frac{t_j^{(N)} - t_k^{(N-1)} + h}
{t_j^{(N)} - t_k^{(N-1)}} 
&  \prod_{k \neq j}
\frac {t_j^{(N)} - t_k^{(N)} - h}
{t_j^{(N)} - t_k^{(N)} + h} \times \\
& \times \prod_{k=1}^{2l_N}\frac{t_j^{(N)} - t_k^{(N)} + h/2}
{t_j^{(N)} - t_k^{(N)} -h/2}=1.
\end{align*}

Note that in the quasiclassical limit $h\to 0$, this system  
becomes system (2.2) of \cite{MV} specialized to the case of 
$C_N$. 



\subsection{Reproduction procedure of $C_1$ type}
Fix a polynomial $T(x)$. 
Suppose that we have a polynomial $y(x)$, such that $(y(x),y(x+h/2))$
is a critical point of $sl_3$ with weights and shifts given by 
$(T(x),T(x+h/2))$ via \Ref{T polyn}. Then there exists a polynomial
$\tilde y$ such that
\be
W(y,\tilde y)(x)=y(x+h/2)T(x).
\ee
The fundamental space $V$ of $sl_3$ population $P_A$ originated at
$(y(x),y(x+h/2)$ is 3-dimensional. It is spanned by $u_1=y(x)$, 
$u_2=\tilde y(x)$ and $u_3$ satisfying the identities
\bean\label{c1 witt}
W(u_1,u_2)(x)=u_1(x+h/2)T(x),\notag\\ 
W(u_1,u_3)(x)=u_2(x+h/2)T(x),\\
W(u_2,u_3)(x)=u_3(x+h/2)T(x).\notag
\eean
These three equations constitute the fact that $u_1, u_2, u_3$ form a
Witt basis of $V$. 

Equations \Ref{c1 witt} imply the following lemma.
\begin{lem}
Let $v$ be a polynomial in $V$. 
The pair $(v(x),v(x+h/2))$ is in the population $P_A$ if
and only if $v(x)$ is a scalar multiple of $u_1+\al u_2+\al^2u_3/2$
for some $\al\in\C$  or  of $u_3$ (i.e., $\al=\infty$). $\square$
\end{lem}
We say that the polynomials $v(x)$
described in the lemma and considered as elements of  $\bs P(\C[x])$  
form a $C_1$ population
of $h$-critical points with weight $T(x)$.
Note that unlike the case of $sl_2$  a $C_1$ population is not a
linear space. 

\subsection{Reproduction procedure of $C_N$ type}
We describe the concepts of reproduction and population. 
All of that follows from the definitions in the cases of
$sl_{2N+1}$ and $C_1$.

An $N$-tuple $\bs y$ is called fertile in $C_N$ sense is there exist 
polynomials
$\tilde y_i$, $i=1,\dots,N-1$ and $\bar y_N$ such that
\bea
W(y_i,\tilde y_i)(x)&=&y_{i+1}(x)y_{i-1}(x+h)T_i(x),\\
W(y_N,\bar y_N)(x)&=&y_{N-1}(x+h)y_N(x+1/2h)T_N(x).
\eea

For $i=1,\dots,N-1$, the
$N$-tuple $\bs y^{(i)}=(y_1,\dots,\tilde y_i,\dots,y_N)$ is called
an immediate descendant of $\bs y$ in the direction $i$.

It follows from the $N$-th equation that the polynomial $y_N$
belongs to a $C_1$ population with weight $y_{N-1}(x+h)T_N(x)$.
Let $\tilde y_N$ be any element of that $C_1$ population. Then
$\bs y^{(N)}=(y_1,\dots, y_{N-1}, \tilde y_N)$ is called
an immediate descendant of $\bs y$ in the direction $N$.       

For $i=1,\dots,N$, immediate descendants $\bs y^{(i)}$
are fertile in the $C_N$ sense.

From Lemma \ref{fertile=critical} we obtain that an $N$-tuple $\bs y$
represents a critical point of $C_N$ type 
if and only if it is fertile and $\bs y^A$
is generic in $sl_{2N+1}$ sense. 
If $\bs y$ represents a critical point of $C_N$ type and 
if $(\bs y^{(i)})^A$ is generic in $sl_{2N+1}$ sense, then $\bs y^{(i)}$ 
also represents an $h$-critical point of type $C_N$.

The $C_N$ population originated at a critical point $\bs y$ is the
minimal set $P(\bs y)$ of fertile $N$-tuples such that $\bs y \in P(\bs y)$ 
and if $\bar{\bs y}\in P(\bs y)$ then the immediate descendants 
${\bar{\bs y}}^{(i)}\in P(\bs y)$. 

Obviously, the $C_N$ population is contained in the corresponding
$sl_{2N+1}$ population: if $\bar{\bs y}\in P(\bs y)$ then  
$\bar{\bs y}^A\in P_A(\bs y^A)$, where we denote $P_A(\bs y^A)$ the
$sl_{2N+1}$ population of critical points originated at $\bs y^A$.
Moreover, the fundamental space of
the $sl_{2N+1}$ population $P_A(\bs y^A)$ contains a flag with the
property \Ref{dar 2 eqn} and therefore is $h$-selfdual. We call this
space the {\it fundamental space of the $C_N$ population $P(\bs y)$}.

The corresponding fundamental operator is the difference operator of
order $2N+1$ obtained from \Ref{Diff} by substituting \Ref{A=C} and
\Ref{T=tilde T}:
\bea
\prod^{N \to 1}_i
\left( \partial  -  
\frac{y_{N-i}(x+(i+3/2)h)}{y_{N-i}(x+(i+1/2)h)} 
\frac{y_{N+1-i}(x+(i-1/2)h)}{y_{N+1-i}(x+(i+1/2)h)}
T(x+ih)\times\right.\\
\left.
\prod_{s=1}^{i-1} \frac{T_{N-s}(x +
  (i-1/2)h)}{T_{N-s}(x+(i-3/2)h)} 
 \right)
\left(\partial -\frac{y_{N}(x+3/2)h)}{y_{N}(x+1/2h)} 
\frac{y_N(x)}{y_N(x+h)}T(x)\right)\times
\\
\prod^{1 \to N}_i
\left( \partial  -  
\frac{y_{N+1-i}(x+h)}{y_{N+1-i}(x)} 
\frac{y_{N-i}(x)}{y_{N-i}(x+h)}
\prod_{s=1}^{N-i} \frac{T_{s}(x + (N-i-s+1)h)}{T_{s}(x+(N-i-s)h)}
 \right),
\eea
where $T(x)=\prod_{s=1}^{N} T_{s}(x + (N-s+1)h)/T_{s}(x+(N-s)h)$.

Clearly two populations of $h$-critical points of type $C_N$ either
do not intersect or coincide. 

\subsection{A $C_N$ population and the $B_N$ flag variety}
Let $\bs y$ be an $h$-critical point of $C_N$ type and let $V$ be the
fundamental space of the population $P_A(\bs y^A)$. Then $V$ is an
$h$-selfdual space of dimension $2N+1$ and therefore has a non-degenerate 
symmetric canonical bilinear form.

The special orthogonal Lie algebra of $V$ consists of all
traceless endomorphisms $x$ of $V$
such that $(xv , v') + (v , xv') = 0$ for all $v, v' \in V$.
Let $\bs u = (u_1, \dots , u_{2N+1})$ be a Witt basis in $V$. We have
$( u_i, u_{2N+2-i}) = (-1)^{i+1}$, $i = 1, \dots , N+1$ ,
and $(u_i , u_j) = 0$ if $i + j \neq 2N+2$.
The choice of the basis identifies the special orthogonal Lie algebra 
with a Lie subalgebra of $sl_{2N+1}$, which is denoted
$so_{2N+1}$. The Lie algebra $so_{2N+1}$ has the root system of type $B_k$.

The lower triangular part of $so_{2N+1}$ is spanned by
matrices $ E_{i, i+1} + E_{ 2N+1-i, 2N+2-i}$  for $i = 1, \dots, N$.
Denote these matrices $X_1, \dots , X_N$, respectively.

Now we describe the action of one-parametric subgroups in the special
orthogonal group in the 
directions of $X_i$ on the basis $\bs u=(u_1,\dots,u_{2N+1})$. 
We have
\bea
&&e^{c X_i} \bs u =\\
&&=(u_1, \dots , u_{i-1}, u_i + c u_{i+1}, u_{i+1}, \dots ,
u_{2N+1-i}, u_{2N+2-i} + c u_{2N+3-i}, \dots , u_{2N+1}),\\
&&e^{c X_i} \bs u =
(u_1, \dots, u_{N-1},u_N + c u_{N+1} + 
c^2 u_{N+2}/2, u_{N+1} + c u_{N+2},
u_{N+2}, \dots , u_{2N+1}),
\eea
where $i=1,\dots,N-1$. 

We also set
\bea
e^{\infty X_i} \bs u &=&
(u_1, \dots , u_{i-1}, u_{i+1}, u_{i}, \dots ,
u_{2N+1-i}, u_{2N+3-i}, u_{2N+2-i}, \dots , u_{2N+1}),\\
e^{\infty X_i} \bs u & = &
(u_1,  \dots , u_{N-1}, u_{N+2}, u_{N+1},u_N
  \dots , u_{2N+1}).
\eea

Note $e^{cX_i}\bs u$ is a Witt basis
for all $i,c$ and the corresponding flag is isotropic.

Denote $FL^\perp(V) \subset FL(V)$ the subvariety of all isotropic
flags and define the $C_N$ generating morphism  
$\beta : FL^\perp(V) \to  \PCN , \mathcal F(\bs u) \mapsto
(y_1^{\mc F},\dots,y_N^{\mc F})$, where 
$y_i^{\mc F}=W^\dagger(u_1,\dots,u_i)$.

\begin{lem}\label{desc lem C}
The set $\beta (\mc F(e^{cX_i}\bs u))$, $c\in\bar{\C}$ 
coincides with the set of all
immediate descendants of $\beta(\mc F(\bs u))$ in the direction $i$.
\end{lem}
\begin{proof}
Follows from the Wronskian identities
\begin{align*}
W(y_i^{\mc
  F},W^\dagger(u_1,\dots,u_{i-1},u_{i+1}))(x)&=T_i(x)y_{i-1}^{\mc
  F}(x+h)y_{i+1}^{\mc F}(x),\\
W(y_N^{\mc F},W^\dagger(u_1,\dots,u_{N-1},u_{N+1}))(x)&=T_N(x)y_{N-1}^{\mc
  F}(x+h)y_N^{\mc F}(x+h/2),\\
W(y_N^{\mc F},W^\dagger(u_1,\dots,u_{N-1},u_{N+2}))(x)&=\\
=T_N(x)y_{N-1}^{\mc F}(x+h)&W^\dagger(u_1,\dots,u_{N-1},u_{N+1})(x+h/2),
\end{align*}
where $i=1,\dots, N-1$.
\end{proof}

\begin{theorem}
The generating morphism $\beta$ of type $C_N$ gives rise 
to the isomorphism $FL^\perp(V)\to P(\bs y)\subset\PCN$.
\end{theorem}
\begin{proof}
The generating morphism $\beta$ is an isomorphism of $FL^\perp(V)$ to the
image by Theorem  \ref{Main}. Clearly, the image of $\beta$ contains the
population $P(\bs y)$. The image of $\beta$ coincides with $P(\bs y)$ by
Lemma \ref{desc lem C}.
\end{proof} 

Note that $FL^\perp$ is isomorphic to the flag variety of the special
symplectic group which corresponds to the $B_N$ root system. This flag
variety has a Bruhat cell decomposition
$FL^\perp(V)=\sqcup_{w\in \mc W}G_w^{\mc F(\infty),B}$, 
where $\mc W$ is the $B_N$ Weyl group.
The Weyl groups of $B_N$ and $C_N$ are canonically identified.
See for details, e.g. \cite{MV}, section 7.4.

As before, the set of all $N$-tuples in a $C_N$ population of the same degree
is the Bruhat cell in the $B_N$ flag variety.
\begin{theorem}
The set of $N$-tuples $\bar{\bs y}$ in a $C_N$ population 
$P(\bs y)$ associated with a
weight at infinity $w\cdot\La_\infty$, where $\La_\infty$ is integral
dominant, coincides with the image of Bruhat cell
$\beta(G^{\mc F(\infty),B}_w)$. 
\end{theorem}
\begin{proof}
Completely parallel to the proof of Corollary 7.14 in \cite{MV}.
\end{proof}

\section{Appendix A: an example of an $sl_3$ population} 
Consider the population of $h$-critical points associated to
$N=2$ and $n=0$ and originated at $\bs y^0 = (1, 1)$. The pair $(1, 1)$
represents the $h$-critical point with no variables. 
We claim that this population consists of pairs of non-zero 
polynomials $\bs y = (y_1, y_2)$, where 
\bean\label{Example}
y_i =  a_{2,i} x^2 +  a_{1,i}x+  a_{0,i} , \qquad i=1, 2 \ ,
\eean
and
\be
(a_{1,1} + a_{2,1} h)  ( a_{1,2} - a_{2,2} h ) 
= 2 a_{0,1}  a_{2,2} +  2 a_{2,1}  a_{0,2} .
\ee
For any  generic pair $\bs y = (y_1, y_2)$ of this form, 
the roots of the polynomials $y_1, y_2$ 
form an $h$-critical point with the initial data where $n=0$.
In other words, the roots of
$y_1, y_2$ satisfy the equations
\bea
\prod_{k \neq j}
\frac {t_j^{(1)} - t_k^{(1)} - h}
{t_j^{(1)} - t_k^{(1)} + h} \
\prod_{k=1}^{l_{2}}
\frac{t_j^{(1)} - t_k^{(2)}}
{t_j^{(1)} - t_k^{(2)} - h} \ = \ 1 \ ,
\qquad j = 1, \dots , l_1 ,
\\
\prod_{k=1}^{l_{1}}
\frac{t_j^{(2)} - t_k^{(1)} + h}
{t_j^{(2)} - t_k^{(1)}} \
\prod_{k \neq j}
\frac {t_j^{(2)} - t_k^{(2)} - h}
{t_j^{(2)} - t_k^{(2)} + h} \ = \ 1\ ,
\qquad j = 1, \dots , l_2 ,
\eea
where $l_1 = $ deg $y_1$ and $l_2 = $ deg $y_2$.

Equations (\ref{main}) take the form
\bean\label{example}
W(y_1, \tilde{y}_1)(x) =  y_2(x) , \qquad W(y_2, \tilde{y}_2)(x) =  y_1(x+h),
\eean
and the reproduction procedure works as follows. We start with $\bs
y^0 = (1, 1)$. 
Equations (\ref{example}) have the form 
$W( 1,\ \tilde{y}_1)(x) =  1$, \ $ W( 1,\ \tilde{y}_2)(x) =   1$.
Using the second of them, we get pairs $\bs y = (1, \ x + a )$ for all
numbers $a$.  
Equations (\ref{example}) now are
$W( 1, \ \tilde{y}_1) =  x + a$, \ $ W( x + a,\ \tilde{y}_2) =  1$.
Using the first equation we get pairs 
$\bs y = ( x^2 + (2a - h)x + b, \ x + a)$ for all  $a, b$.
Equations (\ref{example}) take the form 
$W( x^2 + (2a - h)x + b, \ \tilde{y}_1) = x + a$, \ 
$W(  x + a, \ \tilde{y}_2) =  (x+h)^2 + (2a - h)(x+h) + b$.
Using the second of them  we get 
$\bs y = (x^2 + (2a - h)x + b, \ x^2 + c x + ac - ah - b )$ for all
$a, b, c$.  
It is easy to see that the
union of all those pairs is our population, and nothing else can be
constructed 
starting from $\bs y^0 = (1, 1)$. 

Now, it is easy to see that the family of pairs (\ref{Example}) (where 
each pair is considered up to multiplication of its coordinates by 
non-zero numbers) is isomorphic as an algebraic variety to the variety of all
full flags in the three dimensional vector space $V=\C_2[x]$ of the first
coordinates of the pairs. 
Namely, $y_1$ generates a line in $V$ and $y_2$ defines a plane in 
$V$ containing the line generated by $y_1$.
The space $V$ is the fundamental space of the population.

In the example above the possible degrees of polynomials $y_1, y_2$
are (0,0), (1,0), (0,1), (1,2), (2,1), (2,2). The corresponding parts
of the family are isomorphic to open Bruhat cells of dimensions 0, 1,
1, 2, 2, 3, respectively.

The fundamental space of the population $V$ 
is the space of polynomials of
degree at most 2. This is an $h$-selfdual space. The canonical form in the
basis $(1,x,x(x-1)/2)$ is described by the matrix:
\be
\left(\begin{matrix} 
0 & 0 & 1\\
0 & -1 & 1/2 \\
1 & 1/2 & -1/8
\end{matrix}\right).
\ee
It is a symmetric non-degenerate form. 

The basis $(1,x-1/2,(x^2-x-1/8)/2)$ is a Witt basis. 

The set of isotropic vectors is the $C_1$ population. In this example,
this is the set of polynomials of the form 
\be
(x(x-1)-1/8)/2+\al(x-1/2)+\al^2/2=(x+\al-1/2)^2-1/8,
\ee
where $\al\in \C$ and also $1\in\V$ which corresponds to $\al=\infty$.

\section{Appendix B: the Wronskian identities}
In this appendix we collect identities involving discrete Wronskians.
All functions in this section are functions of one variable $x$.

Set
\be
\Delta f (x) =\ f(x+h) - f(x).
\ee
Set $\Delta^{(0)}f(x) = f(x)$ and
\be
\Delta^{(n+1)} f(x) \ =\ \Delta (\Delta^{(n)} f ) (x)\ .
\ee
The discrete Wronskian of functions $g_1, \dots , g_s$ is defined by
\be
W (g_1, \dots , g_s)(x) = \det(g_i(x + (j - 1)h))_{i,j=1}^s =
\det( \Delta^{(j-1)} g_i(x))_{i,j=1}^s\ .
\ee
We follow the convention that for $s=0$ the corresponding discrete
Wronskian equals $1$. 
In this section we write $W_s(g_1, \dots , g_s)$ instead of $W(g_1, \dots , g_s)$ to
stress the order of the Wronskian.

\begin{lem}\label{1 wronskian}
We have $W_{s+1}(1, g_1, \dots , g_s)(x)\ =
\ W_s(\Delta g_1, \dots , \Delta g_s)(x)$.
\hfill $\square$
\end{lem}

Now we describe the Wronskian 
$W_s(\Delta g_1, \dots , \Delta g_s)(x)$ in more detail. Apriori, this
Wronskian consists of $2^s s!$ terms. However, there are
cancellations. We describe the surviving terms.

\begin{lem}\label{delta wronskian}
The Wronskian 
$W_s(\Delta g_1, \dots , \Delta g_s)(x)$ is equal to 
the alternating sum of 
$(s+1)!$ terms of the form $\prod_{i=1}^s g_i(x+\om_ih)$ 
where
$(\om_1,\dots,\om_s)$ is a permutation $\om$ 
of the set $\{0,\dots,\hat j,\dots, s\}$ for
some $j$. The sign of this term is $(-1)^j\on{sgn} (w)$.
\end{lem}
\begin{proof} 
This lemma is straightforward.
\end{proof}

Now we proceed to other identities.

\begin{lem}\label{f wronskian}
We have $ W_s(fg_1, \dots , fg_s)(x)\ =\  W_s(g_1, \dots , g_s)(x)\
\prod_{j=0}^{s-1} f(x+jh)$ .
\hfill $\square$
\end{lem}

For given functions $g_1, \dots , g_{s+1}$ and an  integer $k$,
$0 \leq k \leq s$, denote
$V_{s-k+1}(i)\ =\ W_{s-k+1}(g_1, \dots , g_{s-k}, g_i)$.

\begin{lem}\label{wr id 2} We have
\bea
&&W_{k+1}(V_{s-k+1}(s-k+1), \dots ,  V_{s-k+1}(s+1))(x) =
\\ 
&&
{\ } \ {}
\qquad
W_{s+1}(g_1, \dots , g_{s+1})(x)\ \prod_{j=1}^k W_{s-k}(g_1, \dots , g_{s-k})(x+jh)\ . 
\eea
\end{lem}
\begin{proof}
This lemma is proved by induction on $s$. The case of $k=s$ is
trivial. Suppose that the lemma is proved
for $s=k,  k+1, \dots, s_0 - 1$. Divide both sides of the 
identity for $s = s_0$
by $\prod_{j=0}^{s_0 - k}\prod_{i=0}^k g_1(x+(i+j)h)$ and 
use Lemma \ref{f wronskian} to carry $g_1$
inside  all Wronskians. Then one of the functions in each Wronskian
is 1 and we can reduce the order by Lemma \ref{1 wronskian}.
Then the identity for $s = s_0$ follows from the  induction 
hypothesis applied to $f_i = \Delta (g_{i+1}/g_1)$, $i = 1, \dots, s_0$.
\end{proof}

For given functions $g_1, \dots , g_{s+1}$, denote
$W_s(i)\ =\ W(g_1, \dots , \widehat{g_i}, \dots , g_{s+1})$
the Wronskian of all functions except  $g_i$.
\begin{lem}\label{wr id 1}
We have
\bea
&&
W_{k+1}(W_s(s+1), W_s(s), \dots ,
W_s(s-k+1))(x)\ =
\\
&&
{\ } \ {\ }
\qquad
 W_{s-k}(g_1, \dots , g_{s-k})(x + kh )\ 
\prod_{j=1}^k W_{s+1}(g_1, \dots , g_{s+1})(x + (j - 1)h)\ .
\eea
\end{lem}
\begin{proof}
  First we prove the case of $s=k$ by induction on $k$. The case of $k=0$ is 
  trivial. Suppose the case of $k<k_0$ is proved. Divide both sides of
  our identity in the case of $s=k=k_0$ by 
$
\prod_{i=0}^{k_0-1}\prod_{j=0}^{k_0}g_1(x+(i+j)h). 
$
By Lemmas
  \ref{1 wronskian} and \ref{f wronskian} we are reduced to the
  identity
\be
W_{k_0+1}(W^{\Delta h}_{k_0-1}[k_0],\dots,W^{\Delta h}_{k_0-1}[1],
W_{k_0}(h_1,\dots,h_{k_0}))
 = \prod_{i=0}^{k_0-1}W(\Delta h_1,\dots,\Delta h_{k_0})(x+ih),
\ee
where $h_i = g_{i+1}/g_1$ and $W^{\Delta h}_{k_0-1}[i] = W_{k_0-1}(\Delta 
h_1,\dots,\widehat{\Delta h_i},\dots,\Delta h_{k_0})$.

The left hand side of the last identity is a determinant of size $k_0+1$. Add to the last
row the row number $i$ with coefficient $(-1)^{k_0-i+1}h_{k_0-i+1}(x+k_0h)$,
$i=1,\dots,k_0$. Then using Lemma \ref{delta wronskian} we observe
that the last row becomes
\be
(0,\dots,0,W_{k_0}(\Delta h_1,\dots,\Delta h_{k_0})(x+(k_0-1)h))
\ee
and the lemma for $k=k_0$ follows from the induction hypothesis
applied to functions $f_1=\Delta h_1,\dots,f_{k_0}=\Delta h_{k_0}$.

Now we continue by induction on $s$. Suppose that the lemma is proved
for $s=k,\dots,s_0-1$. Divide both sides of the identity for $s=s_0$
by $\prod_{i=0}^{s_0-1}\prod_{j=0}^k g_1(x+(i+j)h)$. 
Then the identity for $s=s_0$ follows from the  induction
  hypothesis applied to $f_i=\Delta(g_{i+1}/g_1)$, $i=0,\dots,s_0-1$.
\end{proof}


\begin{thebibliography}{00000}
\normalsize


\bibitem[BIK] {BIK} N.M. Bogoliubov, A.G. Izergin and V.E. Korepin,
  Quantum Inverse Scattering Method and Correlation Functions,
  Cambridge University Press, Cambridge, 1993.

\bibitem[F]{F} W. Fulton, Intersection Theory, Springer-Verlag, 1984.

\bibitem[Fd]{Fd} L.D. Faddeev, {\it Lectures on quantum inverse
    scattering method}, in Integrable Systems, Nankai Lectures on Math
    Phys., 1987 (X-C. Song. ed.), World Scientific, Singapore (1990), 23-70.


\bibitem[FT]{FT} L.D. Faddeev, L.A. Takhtajan , {\it Quantum inverse problem
    method and the Heisenberg $XYZ$ model}, Russian Math. Survey, 34
    (1979), {\bf 5}, 11-68.

\bibitem[GH]{GH} Ph. Griffiths, J. Harris, Principles of
    algebraic geometry, A Whiley-Interscience Publication, 1994.





\bibitem[MV1]{MV} E. Mukhin and A. Varchenko,
{\it Critical Points of Master Functions and Flag Varieties}, 
math.QA/0209017 (2002), 1-49.

\bibitem[MV2]{MV2} E. Mukhin and A. Varchenko,
{\it The quantized Knizhnik-Zamolodchikov equation in tensor products of
irreducible sl(2)-modules}, Calogero-Moser-Sutherland Models workshop,
Springer (2000), 347-384.

\bibitem[TV1]{TV} V. Tarasov and A. Varchenko,
{\it Geometry of q-hypergeometric functions as a bridge between Yangians
and quantum affine algebras}, Invent. math. 128, (1997) 501-588.

\bibitem[TV2]{TV2} V. Tarasov and A. Varchenko,
{\it Completeness of Bethe vectors and Difference equations with 
Regular Singular points}, International Mathematics Research Notices
(1995), {\bf no. 13}, 637-669.






\end{thebibliography}
\end{document}